\documentclass[preprint]{elsarticle} 

\makeatletter
\def\ps@pprintTitle{%
	\let\@oddhead\@empty
	\let\@evenhead\@empty
	\def\@oddfoot{\centerline{\thepage}}%
	\let\@evenfoot\@oddfoot}
\makeatother

\usepackage{etoolbox}
\patchcmd{\MaketitleBox}{\footnotesize\itshape\elsaddress\par\vskip36pt}{\footnotesize\itshape\elsaddress\par\parbox[b][36pt]{\linewidth}{\vfill\hfill\textnormal{\today}\hfill\null\vfill}}{}{}%
\patchcmd{\pprintMaketitle}{\footnotesize\itshape\elsaddress\par\vskip36pt}{\footnotesize\itshape\elsaddress\par\parbox[b][36pt]{\linewidth}{\vfill\hfill\textnormal{\today}\hfill\null\vfill}}{}{}%

\usepackage[
colorlinks=true,%
breaklinks=true,%
linkcolor=blue,
urlcolor=blue,%
citecolor=blue,%
pdftitle={Title of paper}, 
pdfkeywords={Keywords},	
pdfauthor={Authors},		
bookmarksopen=false,
pdfpagemode=None]{hyperref}
\usepackage{algorithm}
\usepackage{algorithmicx}
\usepackage{algpseudocode}
\usepackage{tabularx}
\usepackage{array}
\usepackage{multirow}
\usepackage{threeparttable}
\usepackage{booktabs}
\usepackage{graphicx}
\usepackage{arydshln}

\usepackage[margin = 1.2in]{geometry}
\usepackage[T1]{fontenc}
\usepackage[english]{babel}
\usepackage[latin1]{inputenc}
\usepackage{mathtools}
\usepackage{amsmath}
\usepackage{amsfonts}
\usepackage{amsthm}
\usepackage{amssymb}
\usepackage{array}
\usepackage{subcaption}
\usepackage{siunitx}

\usepackage{color}
\usepackage{url}
\usepackage{cleveref}
\usepackage{graphicx}
\usepackage{breakcites}
\usepackage{soul} 
\usepackage{enumitem}
\usepackage[title,titletoc,toc]{appendix}

\usepackage[textsize=tiny]{todonotes}
 
 \usepackage{setspace}

 \newtheorem{assumption}{Assumption}
{\noindent {\textbf{Proof}:} }%
{\hfill $\Box$ \\[1ex] }

\newcommand{\bit}{\begin{itemize}}
\newcommand{\eit}{\end{itemize}}
\newcommand{\ben}{\begin{enumerate}}
\newcommand{\een}{\end{enumerate}}

\newcommand {\real} {\mathbb{R}}
\newcommand {\nat} {\mathbb{N}}

\DeclareMathOperator*{\Exp}{\mathbb{E}}

\DeclareMathOperator*{\argmin}{arg\,min}%
\newcommand{\rd}{\text{\upshape d}} 

\newcommand{\tk}{\ensuremath{\tilde{k}}}

\newcommand{\bA}{\ensuremath{\mathbf{A}}}

\newcommand{\bG}{\ensuremath{\mathbf{G}}}

\newcommand{\bM}{\ensuremath{\mathbf{M}}}

\newcommand{\bW}{\ensuremath{\mathbf{W}}}
\newcommand{\bX}{\ensuremath{\mathbf{X}}}

\newcommand{\bb}{\ensuremath{\mathbf{b}}}
\newcommand{\bc}{\ensuremath{\mathbf{c}}}

\newcommand{\bj}{\ensuremath{\mathbf{j}}}

\newcommand{\bx}{\ensuremath{\mathbf{x}}}

\newcommand{\cB}{\ensuremath{\mathcal{B}}}

\newcommand{\cD}{\ensuremath{\mathcal{D}}}

\newcommand{\cF}{\ensuremath{\mathcal{F}}}

\newcommand{\cJ}{\ensuremath{\mathcal{J}}}

\newcommand{\cU}{\ensuremath{\mathcal{U}}}

\newtheorem{remark}{Remark}

\newtheorem{problem}{Problem}

\begin{document}
	
	\begin{frontmatter}
		
		\title{Bi-fidelity conditional value-at-risk estimation by dimensionally decomposed generalized polynomial chaos expansion}
		
		
		\author[affil1]{Dongjin Lee\corref{cor1}}
		\ead{dongjin-lee@ucsd.edu}		
		\author[affil1]{Boris Kramer\corref{cor1}}
		\ead{bmkramer@ucsd.edu}
		
		\cortext[cor1]{Corresponding author}
		
		\address[affil1] {Department of Mechanical and Aerospace Engineering, University of California San Diego, CA, United States}

		
		\begin{abstract}
			Digital twin models allow us to continuously assess the possible risk of damage and failure of a complex system. Yet high-fidelity digital twin models can be computationally expensive, making quick-turnaround assessment challenging. 
			Towards this goal, this article proposes a novel bi-fidelity method for estimating the conditional value-at-risk (CVaR) for nonlinear systems subject to dependent and high-dimensional inputs. 
			For models that can be evaluated fast, a method that integrates the dimensionally decomposed generalized polynomial chaos expansion (DD-GPCE) approximation with a standard sampling-based CVaR estimation is proposed.
			For expensive-to-evaluate models, a new bi-fidelity method is proposed that couples the DD-GPCE with a Fourier-polynomial expansion of the mapping between the stochastic low-fidelity and high-fidelity output data to ensure computational efficiency. 
			The method employs measure-consistent orthonormal polynomials in the random variable of the low-fidelity output to approximate the high-fidelity output. 
			Numerical results for a structural mechanics truss with 36-dimensional (dependent random variable) inputs indicate that the DD-GPCE method provides very accurate CVaR estimates that require much lower computational effort than standard GPCE approximations. 
			A second example considers the realistic problem of estimating the risk of damage to a fiber-reinforced composite laminate. The high-fidelity model is a finite element simulation that is prohibitively expensive for risk analysis, such as CVaR computation. Here, the novel bi-fidelity method can accurately estimate CVaR as it includes low-fidelity models in the estimation procedure and uses only a few high-fidelity model evaluations to significantly increase accuracy.
		\end{abstract}	
		
		
		
	\end{frontmatter}
	
\section{Introduction} \label{intro}
Risk assessment is essential for designing and maintaining high-performance engineering systems from the conceptual design stage to operation---where digital twins play a growing role---to product retirement. 
Measures of risk (rather than reliability) have been employed in finance for quite some time, see \citep*{artzner1999coherent,RTRockafellar_SUryasev_2002a,sarykalin2008value}. For example, the Value-at-Risk (VaR) is based on a quantile of the distribution of the output quantity of interest to measure aggregate losses, while the Conditional Value-at-Risk (CVaR) reflects the mean or average size of losses exceeding the VaR. For optimization in either portfolio management \citep*{rockafellar2000optimization,RTRockafellar_SUryasev_2002a,RMansini_WOgryczak_MGSperanza_2007a} or engineering design \citep*{HYang_MGunzburger_2016a,JORoyset_LBonfiglio_GVernengo_SBrizzolara_2017a}, CVaR is superior to VaR in that it quantities tail risk and, as a coherent risk measure, is subadditive \citep*{sarykalin2008value}. Moreover, in contrast to VaR, CVaR preserves convexity of the function it is applied to, which facilitates optimization \citep*{rockafellar2000optimization,kouri2016risk,garreis2021interior,Chaudhuri2022}. Despite offering convexity, a drawback of CVaR is that it is non-smooth, yet smoothed approximations exist~\cite{kouri2016risk,kouri2020epi}, which significantly improve optimization performance. 
It has also been noted that CVaR has several quantitative and qualitative advantages over reliability (or failure probability), a commonly used concept in engineering practice, see \citep*{tyrrell2015engineering,Chaudhuri2022} for a detailed discussion.
Applications of CVaR in the engineering domains, such as the design of civil \citep*{tyrrell2015engineering}, naval \citep*{JORoyset_LBonfiglio_GVernengo_SBrizzolara_2017a}, and aerospace \citep*{HYang_MGunzburger_2016a,JORoyset_LBonfiglio_GVernengo_SBrizzolara_2017a,chaudhuri2020multifidelity} engineering, have appeared. 
Moreover, recent studies imply that CVaR can improve risk management for digital twins \citep{Luca2021digitaltwin,niederer2021scaling,brunton2021data}, encompassing design, procurement, testing, and production.

Since CVaR is a statistical risk measure, its estimation (besides simple cases where an output follows a distribution from a specific parametric family~\citep*{norton2021calculating}) is mostly performed via a sampling method, such as Monte Carlo Simulation (MCS). 
The CVaR and other risk measures are typically associated with a tail of the output distribution. Consequently, a large amount of (at least thousands of) output samples should be obtained to capture the tail risk. However, in most engineering systems, the output (or response) of a system can only be obtained through experiments or computer simulation. The latter is most commonly done through finite element analysis (FEA) or other spatial discretization techniques which can require high computational cost. This makes CVaR estimation computationally intensive if not prohibitive when crude MCS is employed. This situation is compounded when CVaR is used within an optimization problem. 
To mitigate some of these computational challenges, numerous surrogate methods coupled with other variance reduction techniques for CVaR estimation have been developed, such as reduced-order models (ROMs) \citep*{HKTQ18CVaRROMS,ZZou_DPKouri_WAquino_2019a,HKT2020_Adaptive_ROM_CVAR_estimation}, polynomial chaos expansion (PCE) \citep*{bernal2020volatility}, Kriging \citep*{chen2012stochastic}, support vector  machine \citep*{gotoh2017support}, and neural networks \citep*{soma2020statistical}. In addition, recent work \citep*{JAKEMAN2022108280} presents a surrogate modeling method to train the PCE or other surrogate models using limited samples for conservatively estimating CVaR. The method constructs a surrogate model that is tailored to the user's risk preferences (biased to risk measure) while allowing for overestimating risk. 
That work, and most existing other methods, make the simplifying assumption that the input random variables are statistically independent, which then allows factoring their joint probability distribution as the tensor product of the marginal probability distributions of the input variables.   
However, in practice, input variables are often correlated or dependent.  Indeed, neglecting the correlation in input random variables, whether emanating from loads, material properties, or manufacturing variables may produce inaccurate or unknown risky designs \citep*{noh09,lee2020practical}. 

Standard MCS can be used to sample directly from dependent random variables, yet MCS requires only high-fidelity model evaluations and therefore can be computationally prohibitive.
A few other estimation methods, such as generalized PCE (GPCE) \citep*{rahman2018polynomial}, generalized polynomial dimensional decomposition (GPDD) \citep*{rahman19}, or other PCE variants \citep*{Navarro2014POLYNOMIALCE,JAKEMAN2019643}, can handle dependent random variables directly without a potentially detrimental measure transformation between dependent and independent variables. 
A practical version of the GPCE was recently introduced to effectively solve UQ and design optimization problems under arbitrary, dependent input random variables \citep*{lee2020practical,lee2021rdo,lee2021rbdo}. This work makes it possible to obtain the multivariate orthonormal polynomial basis consistent with any non-product-type probability measure of input numerically, instead of an analytical expression by a Rodrigues-type formula used in the prequel \citep*{rahman2018polynomial}. Most recently, a dimensionally decomposed GPCE (DD-GPCE) \citep*{lee2021thesis} has been introduced to tackle stochastic design problems with high-dimensional inputs. As a restructured version of GPCE, the DD-GPCE has been proven to alleviate the curse of dimensionality to some extent by reshuffling and pruning GPCE basis functions in a dimension-wise manner. Two current shortcomings of the DD-GPCE are that, first, the DD-GPCE method has been evaluated only for statistical moment and reliability analyses in design problems. Yet, as elaborated above, in engineering design, risk measures such as CVaR are an interesting alternative and often superior. Second, the DD-GPCE often mandates hundreds of high-fidelity model evaluations which can be computationally prohibitive depending on the complexity and state dimension of the computational models (often FEA models).

This study therefore focuses on CVaR estimation of nonlinear and high-dimensional systems under dependent random variables and enables scalability both with respect to the high state dimension and the high input parameter dimension. Specifically, we propose a novel bi-fidelity method for CVaR estimation of nonlinear systems with high-dimensional, correlated input random variables and/or nonlinear responses. 
The novel method combines (1) the DD-GPCE approximation of a high-dimensional stochastic output function, (2) an innovative method using Fourier-polynomial expansions of the mapping between the stochastic low-fidelity and high-fidelity output data for efficiently calculating the DD-GPCE, and (3) a standard sampling-based CVaR estimation integrated with the DD-GPCE. In contrast to existing bi- or multi-fidelity methods based on an additive and/or multiplicative correction to the low fidelity output \citep*{Kennedy2000,pepper2021local,peherstorfer2018survey}, the proposed bi-fidelity method employs linear or higher-order orthonormal basis functions consistent with the probability measure of the low-fidelity output to approximate the high-fidelity output, thus achieving nearly exponential convergence rate for the output data. Such Fourier-polynomial approximations demand only a handful of high-fidelity output evaluations. The lower-fidelity outputs are determined by DD-GPCE approximations to reduce the computational costs further.

The paper is organized as follows. Section~\ref{sec:2} discusses mathematical notations and preliminaries, including input and output random variables and alternative expressions of CVaR. Also, brief explanations of the GPCE and DD-GPCE methods are provided. Section~\ref{sec:3}  presents a sampling-based CVaR estimation by the DD-GPCE method. Section~\ref{sec:4}  introduces a novel bi-fidelity method for precise and computationally efficient CVaR estimation that requires only a few expensive high-fidelity model evaluations. Numerical results are reported in Section~\ref{sec:5}. Finally, the conclusions are drawn and future directions outlines in Section~\ref{sec:6}.

\section{Background and related methods} \label{sec:2}

This section presents our problem setup and definitions in Section~\ref{sec:2.1}, discusses alternative CVaR definitions in Section~\ref{sec:2.2} and briefly summarizes the GPCE in Section~\ref{sec:2.3} and the DD-GPCE in Section~\ref{sec:2.4}.  

\subsection{Problem setup and definitions} \label{sec:2.1}
Let $\nat$, $\nat_{0}$, $\real$, and $\real_{0}^{+}$ be the sets of positive integers, non-negative integers, real numbers, and non-negative real numbers, respectively. For a positive integer $N\in\nat$, denote by $\mathbb{A}^N \subseteq \real^N$ a bounded or unbounded sub-domain of $\real^N$. 

\subsubsection{Input random variables} \label{sec:2.1:1}  
Let $(\Omega,\cF,\mathbb{P})$ be a probability triple, where $\Omega$ is a sample space representing an abstract set of elementary events, $\cF$ is a $\sigma$-algebra on $\Omega$, and $\mathbb{P}:\cF\to[0,1]$ is a probability measure. Then, consider an $N$-dimensional random vector $\bX:=(X_{1},\ldots,X_{N})^\intercal$, describing the statistical uncertainties in all input and system parameters of a stochastic or random problem. Every so often, $\bX$ will be referred to as an input random vector or input random variables where the integer $N$ represents the total number of input random variables.

Denote by $F_{\bX}({\bx}):=\mathbb{P}\big[\cap_{i=1}^{N}\{ X_i \le x_i \}\big]$ the joint distribution function of $\bX$, admitting the joint probability density function $f_{\bX}({\bx}):={\partial^N F_{\bX}({\bx})}/{\partial x_1 \cdots \partial x_N}$. Given the abstract probability space $(\Omega,\cF,\mathbb{P})$, the image probability space is $(\mathbb{A}^N,\cB^{N},f_{\bX}\rd\bx)$, where $\mathbb{A}^N$ can be viewed as the image of $\Omega$ from the mapping $\bX:\Omega \to \mathbb{A}^N$ and $\cB^N:=\cB(\mathbb{A}^N)$ is the Borel $\sigma$-algebra on $\mathbb{A}^N\subset \mathbb{R}^N$.

We make the following assumptions, which are identical to \citep*{rahman2018polynomial}.  
\begin{assumption} 	\label{a1}
	The random vector $\bX:=(X_{1},\ldots,X_{N})^\intercal$
	\begin{enumerate}
		\item
		has an absolutely continuous joint distribution function $F_{\bX}({\bx})$ and a continuous joint probability density function $f_{\bX}({\bx})$ with a bounded or unbounded support $\mathbb{A}^N \subseteq \real^N$;
		\item
		possesses absolute finite moments of all orders, that is, for all $\bj:=(j_1,\ldots,j_N) \in \nat_0^N$, it holds that
		\begin{equation}
			\Exp \left[ | \bX^{\bj} | \right] :=
			\int_{\Omega}|\bX(\omega)|^{\bj}\mathrm{d}\mathbb{P}(\omega) =
			\int_{\mathbb{A}^N} | \bx^{\bj} | f_{\bX}({\bx}) \rd\bx < \infty,
		\end{equation}
		where $\bX^{\bj}=X_1^{j_1}\cdots X_N^{j_N}$ and $\Exp$ is the expectation operator with respect to the probability measure $\mathbb{P}$ or $f_{\bX}(\bx) \rd\bx$; 
		\item
		has a joint probability density function $f_{\bX}({\bx})$, which
		\begin{enumerate}
			\item
			has a compact support, that is, there exists a compact subset $\mathbb{A}^N \subset \real^N$ such that $\mathbb{P}[\bX \in \mathbb{A}^N]=1$, or
			\item
			is exponentially integrable, that is, there exists a real number $\alpha > 0$ such that
			\begin{equation}
				\int_{\mathbb{A}^N}
				\exp{\left( \alpha \| \bx \|\right)
					f_{\bX}(\bx) \rd\bx} < \infty,
			\end{equation}
			where $\|\cdot\|:\mathbb{A}^N \to \real_0^+$ is an arbitrary norm.
		\end{enumerate}
	\end{enumerate}
\end{assumption}

\subsubsection{Output random variable} \label{sec:2.1:2}
Given an input random vector $\bX:=(X_{1},\ldots,X_{N})^\intercal:\Omega~\to~\mathbb{A}^{N}$
with a known probability density function $f_{\bX}({\bx})$ on $\mathbb{A}^N \subseteq \real^N$, denote by $y(\bX):=y(X_{1},\ldots,X_{N})$ a real-valued, square-integrable output random variable.  Here, $y:\mathbb{A}^N \to \real$ describes a quantity of interest that an application engineer deems relevant for risk assessment. In this work, $y$ is assumed to belong to  the weighted $L^2$ space
\[
\left\{y:\mathbb{A}^N \to \real:~
\int_{\mathbb{A}^N} \left| y(\bx)\right|^2  		f_{\bX}({\bx})\rd\bx < \infty \right\},
\]
which is a Hilbert space. This is tantamount to saying that, for the abstract probability space $(\Omega,\cF,\mathbb{P})$, the output random variable $Y=y(\bX)$ belongs to the equivalent Hilbert space 
\[
L^2(\Omega,\cF,\mathbb{P}):=
\left\{Y:\Omega \to \real:
~\int_{\Omega} \left|y(\bX(\omega))\right|^2 \rd\mathbb{P}(\omega) < \infty\right\}. 
\]
If there is more than one output variable, then each component is associated with a measurement function $y_i$.  Indeed, the generalization for a multivariate output random vector is straightforward.

\subsection {Conditional Value-at-Risk} \label{sec:2.2}
Given a random input $\bX=(X_1,\ldots,X_N)^{\intercal}$, consider an output function $y(\bX)\in L^2(\Omega,\cF,\mathbb{P})$.  For a given risk level $\beta\in(0,1)$, denote by $\mathrm{CVaR}_{\beta}[y(\bX)]$ and $\mathrm{VaR}_{\beta}[y(\bX)]$ the conditional value-at-risk and the value-at-risk of $y(\bX)$ at level $\beta$, respectively. The $\mathrm{VaR}_{\beta}[y(\bX)]$ is the $\beta$-quantile of $y(\bX)$, i.e.,     
\begin{equation}
	\mathrm{VaR}_{\beta}[y(\bX)]=\argmin_{t\in\real}\{\mathbb{P}[y(\bX)\le t]\ge \beta\},
	\label{2.2:1} 
\end{equation} 
where $ \mathbb{P}[y(\bX)\le t]=\int_{\mathbb{A}^N}\mathbb{I}_{\{y(\bx)\le t\}}(\bx)f_{\bX}(\bx)\rd\bx$.
Here, the indicator function is
\[
\mathbb{I}_{\{y(\bx)\leq t\}}(\bx)=\begin{cases}
	1, & y(\bx)\le t, \\
	0, & \text{otherwise}.
\end{cases}
\]  

The $\mathrm{CVaR}_{\beta}[y(\bX)]$ is predicated on the mean value of $y(\bX)$ exceeding $\mathrm{VaR}_{\beta}[y(\bX)]$.
There exist several different equivalent definitions of $\mathrm{CVaR}_{\beta}$. Following \cite*{rockafellar2000optimization,RTRockafellar_SUryasev_2002a}, the $\mathrm{CVaR}_{\beta}$ at level $\beta\in(0,1)$ is 
\begin{align}
	\mathrm{CVaR}_{\beta}[y(\bX)]=\argmin_{t\in\real}\left\{t + \dfrac{1}{1-\beta}\Exp[(y(\bX)-t)_+] \right\},
	\label{2.2:2} 
\end{align}
where $(\cdot)_+=\max(\cdot,0)$ and $\Exp$ is the expectation with respect to $f_{\bX}(\bX)\rd\bx$. The minimum of \eqref{2.2:2} on the interval $\mathrm{VaR}_{\beta}[y(\bX)]\leq t \leq \sup\left\{t:\mathbb{P}[y(\bX)\leq t]\leq \beta \right\}$
is determined by inserting $\mathrm{VaR}_{\beta}[y(\bX)]$ into \eqref{2.2:2}, that is, 
\begin{align} 
	\mathrm{CVaR}_{\beta}[y(\bX)]=\mathrm{VaR}_{\beta}[y(\bX)]
	+\dfrac{1}{1-\beta}\Exp\left[(y(\bX)-\mathrm{VaR}_{\beta}[y(\bX)])_+\right].
	\label{2.2:3} 
\end{align} 
If the cumulative distribution function (CDF) $\mathbb{P}[Y\leq y]$ is continuous at $y=\mathrm{VaR}_{\beta}[y(\bX)]$, the equation \eqref{2.2:3} can be simplified, i.e., 
\begin{equation} 
	\mathrm{CVaR}_{\beta}[y(\bX)]=\dfrac{1}{1-\beta}\Exp[y(\bX)\cdot\mathbb{I}_{\{y(\bX)\ge \mathrm{VaR}_{\beta}[y(\bX)]\}}(\bX)].
	\label{2.2:4} 
\end{equation}  

Having the definition of CVaR and the relevant definitions at hand, we can now formally state the problem that is considered in this paper.

\begin{problem}
	Consider a high-dimensional \textit{dependent} random input vector $\bX \in \mathbb{A}^N$ following an arbitrary probability measure $f_{\bX}(\bx)\rd\bx$, and that satisfies Assumption~\ref{a1}. Moreover, we are given an expensive-to-evaluate output quantity of interest $y: \mathbb{A}^N\mapsto \real$. The goal is to compute the $\mathrm{CVaR}_{\beta}[y(\bX)]$ efficiently.
\end{problem}

\subsection{Generalized polynomial chaos expansion} \label{sec:2.3}
A generalized PCE (GPCE) of a square-integrable random variable $y(\mathbf{X})$ is the expansion of $y(\mathbf{X})$ in terms of an orthonormal polynomial basis in the input variables $\bX$. We briefly review GPCE in this section.  
When $\bX=(X_1,\ldots,X_N)^\intercal$ comprises statistically dependent random variables, the resultant probability measure, in general, is not a product-type, meaning that the joint distribution of $\bX$ cannot be obtained strictly from its marginal distributions.  Consequently, measure-consistent multivariate orthonormal polynomials in $\bx=(x_1,\ldots,x_N)^\intercal$ cannot be built from an $N$-dimensional tensor product of measure-consistent univariate orthonormal polynomials.  In this case, a three-step algorithm based on a whitening transformation of the monomial basis can be used to determine multivariate orthonormal polynomials consistent with an arbitrary, non-product-type probability measure $f_{\bX}(\bx)\rd\bx$ of $\bX$, which will be exploited in the Section~\ref{sec:2.4:1} and \ref{appendix:a}.  

Let $\bj:=(j_1,\ldots,j_N) \in \nat_0^N$ be an $N$-dimensional multi-index.  For a realization $\bx=(x_1,\ldots,x_N)^\intercal \in \mathbb{A}^N \subseteq \real^N$ of $\bX$, a monomial in the real variables $x_1,\ldots,x_N$ is the product $\bx^{\bj}=x_1^{j_1}\ldots x_N^{j_N}$ with a total degree $|\bj|=j_1+\cdots+j_N$.  Consider for each $m \in \nat_0$ the elements of the multi-index set 
\[
\cJ_m:=\{ \bj \in \nat_0^N: |\bj|\le m\},
\]
which is arranged as $\bj^{(1)},\ldots,\bj^{(L_{N,m})}$, $\bj^{(1)}=\textbf{0}$, according to a monomial order of choice.  The set $\cJ_m$ has cardinality $L_{N,m}$ obtained as
\begin{align}
	L_{N,m}:=|\cJ_m|=\sum_{l=0}^m \binom{N+l-1}{l}=\binom{N+m}{m}.
	\label{2.3:1}
\end{align}
Let us denote by
\begin{align}
	{\mathbf{\Psi}}_{m}(\bx)=({\Psi}_1(\bx),\ldots,{\Psi}_{L_{N,m}}
	(\bx))^{\intercal}
\end{align}
an $L_{N,m}$-dimensional vector of multivariate orthonormal polynomials that are consistent with the probability measure $f_{\bX}(\bx)\rd{\bx}$ of $\bX$. Consequently, any output random variable $y(\bX)\in L^2(\Omega, \cF, \mathbb{P})$ can be approximated by the $m$th-order GPCE\footnote{The GPCE in \eqref{2.3:2} should not be confused with that of \citep*{xiu02}. The GPCE, presented here, is meant for an arbitrary dependent probability distribution of random input. In contrast, the existing PCE, whether classical \citep*{wiener38} or generalized \citep*{xiu02}, still requires independent random inputs.}
\begin{align}
	y_m(\bX)=\displaystyle\sum_{i=1}^{L_{N,m}}c_i\Psi_i(\bX)
	\label{2.3:2}
\end{align}
of $y(\bX)$, comprising $L_{N,m}$ basis functions with expansion coefficients 
\begin{align}
	c_i:=\displaystyle\int_{\mathbb{A}^N}y(\bx)\Psi_i(\bx)f_{\bX}(\bx)\rd\bx, i=1,\ldots,L_{N,m}.
\end{align} 
Here, the orthonormal polynomials $\Psi_i(\bX)$, $i=1,\ldots,L_{N,m}$, are determined by the three steps in \ref{appendix:a}.  We refer to~\citep*{lee2020practical} for more details.		
The GPCE is referred to as regular GPCE to distinguish it from the DD-GPCE which is introduced next.

\subsection{Dimensionally decomposed generalized polynomial chaos expansion}  \label{sec:2.4}
For problems with high-dimensional inputs (say, $N\geq 20$), the regular GPCE approximation in~\eqref{2.3:2} requires a relatively large number of basis functions due to the growth of $L_{N,m}$ in~\eqref{2.3:1}, which reflects the curse of dimensionality. For example, for a total degree of $m=3$, consider an increase of $N$ from $20$ to $50$. The respective number of the regular GPCE's basis functions exponentially increases from $1,771$ to $23,426$. However, in many real-world applications,  high-variate interaction effects among input variables are often negligible to the output function value of interest~\citep*{rabitz1999,rahman2008}. In such cases, we can leverage the DD-GPCE method to reorder the basis functions of the regular GPCE in a dimension-wise manner that then allows for effectively truncating them to tackle high-dimensional problems. In the next section, the DD-GPCE is briefly summarized, see~\citep*{lee2021thesis} for details.  

\subsubsection{Measure-consistent orthonormal polynomials}  \label{sec:2.4:1}
The DD-GPCE has the ability to effectively select a subset of the basis functions of the regular GPCE based on the degree of interaction among input variables.  Consequently, the method can capture complex nonlinear behavior of the output functions while reducing the exponential growth of the basis functions. The chosen multivariate orthonormal polynomials that are consistent with an arbitrary, non-product-type probability measure $f_{\bX}(\bx)\rd\bx$ of $\bx$ are determined by the three-step process based on a whitening transformation of the monomial basis as follows.

For $N\in\nat$, denote by $\{1,\ldots,N\}$ an index set and $\cU\subseteq\{1,\ldots,N\}$ a subset (including the empty set $\emptyset$) with cardinality $0\leq|\cU|\leq N$.  The complementary subset of $\cU$ is denoted by $\cU^{c}:=\{1,\ldots,N\}\backslash \cU$.  For each $m\in\nat_0$ and $0\leq S \leq N$, we define the reduced multi-index set
\begin{equation} 
	\begin{array}{rcl}
		\cJ_{S,m}:=\left\{\bj=(\bj_{\cU},\textbf{0}_{\cU^c})\in\nat_0^N:\bj_{\cU}\in\nat^{|\cU|},~|\cU|\leq|\bj_{\cU}|\leq m,\right.\left.0\leq |\cU| \leq S\right\},
	\end{array} 
\end{equation} 
which is arranged as $\bj^{(1)},\ldots,\bj^{(L_{N,S,m})}$,  $\bj^{(1)}=\textbf{0}$, according to a monomial order of choice and where $|\bj_{\cU}| := j_{i_1}+\cdots+j_{i_{|\cU|}}$.  Here, $(\bj_{\cU},\textbf{0}_{\cU^c})$ denotes an $N$-dimensional multi-index whose $i$th component is $j_i$ if $i\in \cU$ and $0$ if $i\notin \cU$.
It is straightforward to show that $\cJ_{S,m}$ has cardinality
\begin{equation}
	L_{N,S,m}:=|\cJ_{S,m}|=1+\sum_{s=1}^{S}\binom{N}{s}\binom{m}{s}.
	\label{sr2} 
\end{equation}
The set $\cJ_{S,m}$ represents a subset of $\cJ_m$ determined from the chosen $S$, where only at most $S$-variate basis functions are preserved, that are relevant for the $S$th-variate DD-GPCE approximation. For example, univariate ($S=1$) and bivariate ($S=2$) DD-GPCE approximations demand $Nm+1$ and $N(N-1)m(m-1)/4+Nm+1$ basis functions, respectively, according to \eqref{sr2}. The expression ``$S$-variate'' used for the DD-GPCE approximation means that at most $S$-degree interaction of input variables are included. As an example, univariate and bivariate input variables are $x_i$ and $x_{i_1}x_{i_2}$, where $i_1\neq i_2$, respectively. Thus, the DD-GPCE approximation, the sum of at most $S$-variate basis functions, is strictly an $N$-variate function, which will be explained in detail in the following subsection.  As a result, we have that
\[
L_{N,S,m} \leq L_{N,m},
\]
i.e., the DD-GPCE never has more terms than the standard GPCE; in most cases, it will have significantly less terms.
For $\bx=(x_1,\ldots,x_N)^{\intercal}\in\mathbb{A}^N\subseteq\real^N$ we then define the basis vector for the DD-GPCE as
\[
\mathbf{\Psi}_{S,m}(\bx):=(\Psi_i(\bx),\ldots,\Psi_{L_{N,S,m}}(\bx))^{\intercal},
\]
which is an $L_{N,S,m}$-dimensional vector of multivariate orthonormal polynomials that is consistent with the probability measure $f_{\bX}(\bx)\rd\bx$ of $\bx$. The orthonormal polynomials are determined by the following three steps.  

\begin{enumerate}
	[labelwidth=1.4cm,labelindent=5pt,leftmargin=1.6cm,label=\bfseries Step  \arabic*.,align=left]
	\setcounter{enumi}{0}
	\item Given $0\leq S \leq N$ and $S \leq m < \infty$, create an $L_{N,S,m}$-dimensional column vector
	\begin{equation}	 
		\bM_{S,m}(\bx)=(\bx^{\bj^{(1)}},\ldots,\bx^{\bj^{(L_{N,S,m})}})^{\intercal}
		\label{2.4:1}
	\end{equation}
	
	of monomials whose elements are the monomials $\bx^{\bj}$ for $\bj\in\cJ_{S,m}$ arranged in the aforementioned order.  For $\cU\subseteq \{1,\ldots,N\}$, let $\bx_{\cU}:=(x_{i_1},\ldots,x_{i_{|\cU|}})^{\intercal}$, $1 \leq i_1 < \cdots < i_{|\cU|} \leq N$, be a subvector of $\bx$.  The complementary subvector is defined by $\bx_{\cU^c}:=\bx_{\{1,\ldots,N\}\backslash \cU}$.  Then, for $\bj\in\cJ_{S,m}$,
	\[
	\bx^{\bj}=\bx_{\cU}^{\bj_{\cU}}{\textbf{0}_{\cU^c}}^{\bj_{\cU^c}}=\bx_{\cU}^{\bj_{\cU}}.
	\]  
	Hence, $\bM_{S,m}(\bx)$ is the monomial vector in $\bx_{\cU}=(x_{i_1},\ldots,x_{i_{|\cU|}})^{\intercal}$ of degree $0\leq |\cU| \leq S$ and $|\cU| \leq |\bj_{\cU}| \leq m$. 
	\item
	Construct an $L_{N,S,m} \times L_{N,S,m}$ monomial moment matrix of $\bM_{S,m}(\bX)$, defined as
	\begin{equation}
		\begin{split}
			{\bG}_{S,m}:= \Exp[\bM_{S,m}(\bX)\bM_{S,m}^{\intercal}(\bX)]=
			\int_{\mathbb{A}^{N}}\bM_{S,m}(\bx)\bM_{S,m}^{\intercal}(\bx)
			f_{\bX}(\bx)\rd\bx.
		\end{split}
		\label{2.4:2}
	\end{equation}		
	For an arbitrary PDF $f_{\bX}(\bx)$, the matrix $\bG_{S,m}$ cannot be determined exactly, yet it can be accurately estimated with numerical integration and/or sampling methods~\citep*{lee2020practical}.
	\item
	Select the $L_{N,S,m} \times L_{N,S,m}$ whitening matrix ${\bW}_{S,m}$ from the Cholesky decomposition of the monomial moment matrix ${\bG}_{S,m}$~\citep*{rahman2018polynomial}, leading to
	\begin{equation}
		{\bW}_{S,m}^{-1}{\bW}_{S,m}^{-\intercal}={\bG}_{S,m}.
		\label{2.4:3}
	\end{equation}	
	Then employ the whitening transformation to generate multivariate orthonormal polynomials from
	\begin{equation}
		{\mathbf{\Psi}}_{S,m}(\bx)={\bW}_{S,m} \bM_{S,m}(\bx).
	\end{equation}	
\end{enumerate}
It is straightforward to show that $\bG_{S,m}$ is symmetric and positive definite. 
However, the effectiveness of the three-step algorithm is dependent on reliable construction of a well-conditioned monomial moment matrix. That is, numerical issues with the Cholesky factorization in \eqref{2.4:3} can occur if the estimated matrix has a large condition number.

For the $i$th element $\Psi_i(\bX)$ of the orthonormal polynomial vector $\mathbf{\Psi}_{S,m}(\bX)=
(\Psi_1(\bX),\ldots,\Psi_{L_{N,S,m}}(\bX))^\intercal$, the first- and second-order moments are \citep*{lee2020practical}
\begin{equation}
	\Exp\left[{\Psi}_i(\bX)\right] =
	\begin{cases}
		1, & \text{if}~{i=1}, \\
		0, & \text{if}~{i\ne 1},
	\end{cases}
\end{equation}
and
\begin{equation}
	\Exp\left[{\Psi}_i(\bX){\Psi}_j(\bX)\right] =
	\begin{cases}
		1, & i=j, \\
		0, & i \ne j,
	\end{cases}
\end{equation}
respectively.  These properties are essential to DD-GPCE which we exploit in the next section.
Interested readers may consult Section 6.2.1.4 of \citep*{lee2021thesis} for an illustrative example of measure-consistent orthonormal polynomials and a comparison between the DD-GPCE and the regular GPCE.

\subsubsection{DD-GPCE approximation} \label{sec:2.4:2}
The $S$-variate, $m$th-order DD-GPCE approximation of $y(\bX)$ is
\begin{equation}
	\displaystyle
	y_{S,m}(\bX) =
	\sum_{i=1}^{L_{N,S,m}}
	c_i \Psi_i(\bX)\simeq y
	\label{dd-gpce}
\end{equation}
with expansion coefficients 
\begin{equation}
	\begin{split}
		c_i &:= \int_{\mathbb{A}^N} y(\bx)\Psi_i(\bx) f_{\bX}(\bx) \rd\bx,\quad i=1,\ldots,L_{N,S,m}\ .		
	\end{split}
\end{equation}
The truncation parameters $S$ and $m$ should satisfy $1 \leq S \leq N$ and $S\leq m < \infty$. The DD-GPCE approximation has the property that its basis functions retain the degree of interaction among input variables less than or equal to $S$ and preserves polynomial orders less than or equal to $m$.  

Since the regular GPCE of $y(\bX)\in L^2(\Omega,\cal{F},\mathbb{P})$ converges as $m\rightarrow \infty$ in the mean-square sense---both in probability and in distribution, see~\citep*[Theorem 1]{lee2020practical}---the DD-GPCE $y_{S,m}(\bX)$ converges to $y(\bX)$ in the same sense as the regular GPCE as $S\rightarrow N$ and $m \rightarrow \infty$. When $y$ itself is a polynomial function, the DD-GPCE with the same degree ($S$) and same order ($m$) of $y$ represents $y$ exactly.

\subsubsection{Expansion coefficients} \label{sec:2.4:3}
The expansion coefficients $c_i$, $i=1,\ldots,L_{N,S,m}$, of an $S$-variate, $m$th-order DD-GPCE approximation $y_{S,m}(\bX)$ are determined via high-dimensional integration. For an arbitrary function $y$ and an arbitrary probability distribution of the random input $\bX$, evaluating the integrals exactly from the definition is impossible. Thus, we need to integrate numerically; yet numerical integration via, e.g., a multivariate, tensor-product Gauss-type quadrature rule becomes computationally intensive, if not prohibitive, when $N\geq 10$ (say). For example, if an expansion coefficient is estimated by $N=10$ dimensional numerical integration with a $3$-point rule for each variable, the total number of integrand evaluations is $3^{10}=59,049$.
To address this difficulty, standard least-squares (SLS), as briefly outlined next, can be employed to estimate the coefficients. 

Given the known distribution of a random input $\bX$ and an output function $y:\mathbb{A}^N \to \real$, consider an input-output data set $\{ \bx^{(l)}, y(\bx^{(l)}) \}_{l=1}^{L^{\prime}}$ of size $L^{\prime}\in \nat$  generated by evaluating the quantity of interest $y$ at each input data $\bx^{(l)}$. Various sampling methods such as standard MCS, Quasi-MCS (QMCS), and Latin hypercube sampling (LHS) or various optimal design of experiments can be employed to build the data set such that the observed distribution underlying the data is consistent with the input distribution, ensuring unbiased estimates of the output quantities. The input-output data set is sometimes referred to as the experimental design data, since the inputs are following the input distribution, i.e., $\bx \sim f_{\bX}(\bx)$.  
Using the data set, we can obtain approximate DD-GPCE coefficients $\tilde{\bc} = (\tilde{c}_1,\ldots,\tilde{c}_{L_{N,S,m}})^{\intercal}$ by finding the least-squares solution, i.e., 
\begin{align}
	\tilde{\bc} = \argmin_{\bc \in \real^{L_{N,S,m}}} \Vert \bb - \bA\bc \Vert_2 ,
	\label{sls}
\end{align}
where
\begin{equation}
	\begin{split}
		\bA :=
		\begin{bmatrix}
			\tilde{\Psi}_1(\bx^{(1)}) & \cdots &  \tilde{\Psi}_{L_{N,S,m}}(\bx^{(1)}) \\
			\vdots                   & \ddots &  \vdots                           \\
			\tilde{\Psi}_1(\bx^{(L^{\prime})}) & \cdots &  \tilde{\Psi}_{L_{N,S,m}}(\bx^{(L^{\prime})})
		\end{bmatrix}\ \ \text{and} 
		~~~\bb := (  y(\bx^{(1)}),\ldots, y(\bx^{(L^{\prime})}))^\intercal.\\
	\end{split}
	\label{2.4:4}
\end{equation}
From \eqref{2.4:4}, the terms $\tilde{\Psi}_i(\bx^{(l)})$ are approximations of $\Psi_i(\bx^{(l)})$ resulting from the construction of the monomial moment matrix in Section~\ref{sec:2.4:1}.  According to SLS, the optimal expansion coefficients are obtained by minimizing the mean-squared residual
\begin{equation}
	\hat{e}_{S,m} :=
	\displaystyle
	\frac{1}{L^{\prime}}
	\displaystyle
	\sum_{l=1}^{L^{\prime}}
	\left[
	y(\bx^{(l)}) -
	\displaystyle
	\sum_{i=1}^{L_{N,S,m}} \tilde{c}_i \tilde{\Psi}_i(\bx^{(l)})
	\right]^2.
\end{equation}
The least-squares solution $\tilde{\bc}$ is obtained from the normal equations $\bA^\intercal \bA \tilde{\bc} = \bA^\intercal \bb$, where  the $L_{N,S,m} \times L_{N,S,m}$ matrix $\bA^\intercal \bA$ is referred to as the information or data matrix. If $\bA^\intercal \bA$ is positive-definite then the solution the least-squares approximation of the DD-GPCE coefficients is
\begin{equation}
	\tilde{\bc} = (\bA^\intercal \bA)^{-1} \bA^\intercal \bb.
\end{equation}
When using SLS, the number of experimental data must be larger than the number of coefficients, that is, $L^{\prime}>L_{N,S,m}$. Even if this condition is met, the experimental design must be carefully chosen to ensure that the resulting matrix $\bA^\intercal \bA$ is well-conditioned.

In this work, the data $\bA,\bb$ for the least-squares solution $\tilde{\bc}$ in \eqref{sls} is obtained from MCS. This satisfies the required accuracy of the estimates of the DD-GPCE coefficients in all examples we considered. Alternatively, one may perform several optimal design of experiments which have proven to be effective in the stability of the least-squares solution, see~\citep*{hadigol2018,luthen2021}.

\section{Sampling-based CVaR estimation by dimensionally decomposed GPCE} \label{sec:3} 

This section shows how the DD-GPCE can be leveraged for sampling-based CVaR estimation. We include a discussion of the difference in convergence rate when the two---theoretically equivalent---CVaR expressions are used in practice. We finish the section by presenting a complete algorithm. 
The estimation of $\mathrm{VaR}$ in \eqref{2.2:1} and $\mathrm{CVaR}$ in either \eqref{2.2:3} or \eqref{2.2:4}  for nontrivial examples requires a sampling method, such as MCS, Quasi MCS (QMCS), or Latin hypercube sampling (LHS). We follow the sampling-based MC algorithm from~\citep*[Alg. 2.1]{HKTQ18CVaRROMS}.

The $S$-variate, $m$th-order DD-GPCE $y_{S,m}(\bX)$, presented in Section \ref{sec:2.4},  can be employed as an inexpensive surrogate model that replaces an expensive-to-evaluate function $y(\bX)$. 
Thus, the sampling-based estimation is performed with $y_{S,m}(\bx)$ and we denote the estimates by $\widehat{\mathrm{VaR}}_{\beta}[y_{S,m}(\bX)]$ and $\widehat{\mathrm{CVaR}}_{\beta}[y_{S,m}(\bX)]$. Algorithm~\ref{a1} summarizes all steps of the estimation process.

We employ Algorithm~\ref{al1} for standard MC sampling, where the probability $p^{(l)}=1/L$ for $l=1,\ldots,L$. We note that the algorithm can also be used in the context of importance sampling, in which case $p^{(l)}=\omega(\bx^{(l)})/M$, $l=1,\ldots,M \ll L$, with $\omega(\cdot)$ being the weight function, see~\citep*{HKTQ18CVaRROMS}.   

While the computation of the value at-risk is straightforward (see Steps 3-4 in Algorithm~\ref{a1}), the computation of the CVaR estimate in Step~5 requires further discussion, as outlined next. 
Based on the two equivalent definitions of ${\mathrm{CVaR}}_{\beta}[y(\bX)]$ in \eqref{2.2:3} and \eqref{2.2:4}, one can obtain the estimate 
\begin{align} 
	\begin{array}{rl}
		\widehat{\mathrm{CVaR}}_{\beta}[y_{S,m}(\bX)]=\widehat{\mathrm{VaR}}_{\beta}[y_{S,m}(\bX)]+\dfrac{1}{1-\beta}\dfrac{1}{L}\displaystyle\sum_{l=1}^L\left(y_{S,m}(\bx^{(l)})-\widehat{\mathrm{VaR}}_{\beta}[y_{S,m}(\bX)]  \right)_+
		\label{3.2}
	\end{array} 	
\end{align} 
from \eqref{2.2:3} or one can compute it from  \eqref{2.2:4} as
\begin{align} 
	\begin{array}{rl} 
		\widehat{\mathrm{CVaR}}_{\beta}[y_{S,m}(\bX)]=\dfrac{1}{1-\beta}\dfrac{1}{L}\displaystyle\sum_{l=1}^L y_{S,m}(\bx^{(l)})\mathbb{I}_{\left\{y_{S,m} \ge \widehat{\mathrm{VaR}}_{\beta}[y_{S,m}(\bX)]\right\}}(y_{S,m}(\bx^{(l)})).
		\label{3.3}
	\end{array} 
\end{align}
For the sampling-based estimation employing Algorithm~\ref{al1}, we found that the convergence rate of these two CVaR estimates can be significantly different. The estimate \eqref{3.2} is generally converging faster than the one \eqref{3.3}, so we use \eqref{3.2} in this study. 
\begin{remark}
	Algorithm~\ref{al1} produces unbiased estimates for ${\rm{VaR}}_{\beta}[y(\mathbf{X})]$ and ${\rm{CVaR}}_{\beta}[y(\mathbf{X})]$ from samples $y(\bx^{(l)})$, $l=1,\ldots,L$. For asymptotic convergence properties of the sampling-based estimators $\widehat{\rm{VaR}_{\beta}}[y(\bX)]$ and $\widehat{\rm{CVaR}_{\beta}}[y(\bX)]$ c.f. \cite[Theorem 2.1]{hong2014monte}; for instance the estimates have distinct convergence rates, ${\cal O}(L^{-3/4}(\log L)^{3/4})$ and ${\cal O}(L^{-1}\log L)$, respectively. 
	We note that standard MCS requires $L$ evaluations of the output $y(\bx^{(l)})$ for $l=1,\cdots,L$, which can create prohibitive computational demands and one may not have enough ressources to obtain converged ${\rm VaR}_{\beta}$ and ${\rm CVaR}_{\beta}$ estimates. In contrast, the MCS with the DD-GPCE approximation requires evaluations of simple polynomial functions to obtain $y_{S,m}$. Taking a large number of samples from the DD-GPCE allows us to obtain converged CVaR solutions.
	
	Moreover, the DD-GPCE approximation $y_{S,m}$ of $y$ converges to $y$ as $S\rightarrow N$ and $m \rightarrow \infty$ in the mean-square sense. As a surrogate model, however, the DD-GPCE usually has a model error, which means that $y_{S,m} \approx y$ (except we get equality when $y$ is a polynomial function and the same degree ($S$) and order ($m$) of $y$ are selected for the DD-GPCE) and so, in general, $\widehat{\rm{CVAR}}_{\beta}[y_{S,m}(\bX)]$ does not converge to ${\rm{CVaR}_{\beta}}[y(\bX)]$. Our numerical results in Section~\ref{sec:5} show that the bias is minimal. 
\end{remark} 
\begin{algorithm}
	\caption{Sampling-based estimation of $\mathrm{VaR}_{\beta}$
		and $\mathrm{CVaR}_{\beta}$ by the $S$-variate, $m$th-order DD-GPCE approximation $y_{S,m}$.}
	\label{al1}
		\begin{algorithmic}[1] 
			\Require Samples $\bx^{(l)}=(x_1^{(l)},\ldots,x_N^{(l)})^{\intercal}$, $l=1,\ldots,L\gg 1$, via MCS, QMCS, or LHS with corresponding probabilities  $p^{(l)}=f_{\bX}(\bx^{(l)})\rd\bx^{(l)}$; Set a risk level $\beta \in (0,1)$.
			\Ensure  Estimates $\widehat{\mathrm{VaR}}_{\beta}[y_{S,m}(\bX)]$ and $\widehat{\mathrm{CVaR}}_{\beta}[y_{S,m}(\bX)]$.
			\State Create output samples $y_{S,m}(\bx^{(l)})\simeq y(\bx^{(l)})$, $l=1,\ldots,L$.
			\State Sort values of $y$ in descending order and relabel the samples so that 
			$$
			y_{S,m}(\bx^{(1)})>y_{S,m}(\bx^{(2)})>\cdots>y_{S,m}(\bx^{(L)}),
			$$
			and reorder the probabilities accordingly (so that $p^{(l)}$ corresponds to $\bx^{(l)}$).
			\State Compute the index $\tk_{\beta}\in\nat$ such that 
			$$
			\displaystyle\sum_{l=1}^{\tk_{\beta}-1}p^{(l)} \le 1-\beta < \sum_{l=1}^{\tk_{\beta}}p^{(l)}.
			$$
			\State Set $\widehat{\mathrm{VaR}}_{\beta}[y_{S,m}(\bX)]={y_{S,m}}(\bx^{({\tilde{k}}_{\beta})})$.
			\State Set $\mathrm{CVaR}_{\beta}[y(\bX)]\approx \widehat{\mathrm{CVaR}}_{\beta}[y_{S,m}(\bX)]$ in \eqref{3.2}.  
		\end{algorithmic}
\end{algorithm}

The DD-GPCE relies on SLS or its variants, where the output sample size is usually determined as at least three or four times of the number of basis functions or expansion coefficients. Thereby, for risk measures of high-dimensional inputs (say, $N\ge20$), the DD-GPCE mandates obtaining hundreds of output samples, which can be computationally intensive when that requires expensive high-fidelity model evaluations.  For such cases, we propose a novel bi-fidelity method to efficiently compute the DD-GPCE, as introduced in the following section.    

\section{Bi-fidelity method for CVaR estimation} \label{sec:4}

In practical applications, the output $y(\bX)\in L^2(\Omega,\cF,\mathbb{P})$ requires the simulation of a computational model, e.g., via FEA. This allows the user to choose the level of fidelity. Computationally expensive high-fidelity models produce accurate solutions whereas faster lower-fidelity models, by definition, introduce output bias (or error). Multifidelity methods combine models with multiple fidelities to solve UQ problems and can produce excellent results with provable guarantees at much lower cost, see the survey~\citep*{peherstorfer2018survey} and references therein. 

This section introduces a novel bi-fidelity method that combines the benefits of lower- and high-fidelity models for the construction of DD-GPCE approximations and their subsequent deployment for precise and efficient CVaR estimation under arbitrary dependent random inputs. 
We define relations of high and low-fidelity output in Section~\ref{sec:4.1}, present a novel Fourier polynomial expansion for high-fidelity output in Section~\ref{sec:4.2}, and describe the complete algorithm of the proposed bi-fidelity method for CVaR estimation and its cost in Section~\ref{sec:4.3}. 

\subsection{Relations of high and low-fidelity output} \label{sec:4.1}
For $y(\bX)\in L^2(\Omega,\cF,\mathbb{P})$, we denote by $y_H(\bX)$ and $y_L(\bX)$ the random output estimated by the high- and low-fidelity models of $y$, respectively, which are both functions of the same random input $\bX\in\mathbb{A}^N\subseteq \real^N$ and thus share the identical sample space $\Omega$. 
Let $Y_L:=y_L(\bX)$ whose PDF $f_{Y_L}(y_L)$ is on the domain of $Y_L$, denoted by $\bar{\mathbb{A}}\subseteq \mathbb{R}$.  

Consider a function $h:\bar{\mathbb{A}}\rightarrow\real$ of $Y_L$ that approximates the high-fidelity output $y_H(\bX)$, i.e., 
\begin{align}
	y_H(\bX) \approx h(Y_L).
	\label{4.1}
\end{align}
The mapping $h$ suggests that the relationship between input and output is much simpler than the high-fidelity model $y_H$. 
Thus, we approximate the high-fidelity output via the mapping $h(Y_L)$; we do so with a low-degree Fourier-polynomial expansion with measure-consistent orthonormal polynomials in $Y_L$, which can be shown to have nearly exponential convergence~\citep*{rahman2018polynomial}. 
The Fourier-polynomial expansion will be introduced in the following subsection. 


		%
				%
				
				\subsection{Fourier polynomial expansion to approximate the high-fidelity output} \label{sec:4.2}
				Given a random variable  $Y_{{L}}$ satisfying Assumption~\ref{a1}, the mapping $h(Y_L)$ in \eqref{4.1} can be obtained via its $\bar{m}$th-order Fourier polynomial expansion as
				\begin{align}
					h_{\bar{m}}(Y_{L})=\sum_{i=1}^{\bar{m}+1}b_{i}\Psi_i(Y_{L}),
					\label{4.1:1} 
				\end{align}
				with its expansion coefficients
				\[
				b_i=\displaystyle\int_{\bar{\mathbb{A}}}h(y_L)\Psi_i(y_L)f_{Y_L}(y_L)\rd y_L,~i=1,\ldots,\bar{m}+1.
				\]
				Here, $\Psi_i(Y_L),~i=1,\ldots,\bar{m}+1$ are orthonormal polynomials that are consistent with the probability measure $f_{Y_L}(y_L)\rd y_L$ of the low-fidelity output random variable $Y_L$. These orthogonal polynomials are determined by a three-step process, as described next. 
				
				\subsubsection{Orthonormal polynomial construction} \label{sec:4.2:1}
				For $j\in\nat_0$ and $y_L\in \bar{\mathbb{A}}\subseteq \real$, a monomial in the real variable $y_L$ is $y_L^{j}$ and has a degree $j$. Consider for each $\bar{m}\in\nat_0$ the elements of the ordered index set $\{ j \in \nat_0: j\le \bar{m}\}$. The set has cardinality $\bar{m}+1$. We denote by
				\begin{align}
					{\mathbf{\Psi}}_{\bar{m}}(y_L)=({\Psi}_1(y_L),\ldots,{\Psi}_{\bar{m}+1}(y_L))^{\intercal},
				\end{align}
				the $\bar{m}+1$-dimensional vector of orthonormal polynomials that are consistent with the probability measure $f_{Y_L}(y_L)\rd{y_L}$ of $Y_L$.  This vector is determined as follows.
				
				\begin{enumerate}
					[labelwidth=1.4cm,labelindent=5pt,leftmargin=1.6cm,label=\bfseries Step  \arabic*.,align=left]
					\setcounter{enumi}{0}
					\item
					Given $\bar{m} \in \nat_0$, $m \in \nat_0$, and $0\leq S \leq m$, determine an $(\bar{m}+1)$-dimensional column vector
					\begin{equation}	 \bM_{\bar{m}}(y_L)=(1,{y_L},y_L^2,\ldots,y_L^{\bar{m}})^{\intercal},
					\end{equation}
					of monomials $y_L^{j}$ for $|j|\le \bar{m}$. The real variable $y_L$ can be replaced with the $S$-variate, $m$th-order DD-GPCE   
					\begin{align}
						\tilde{y}_{L,S,m}(\bx)=\sum_{i=1}^{L_{N,S,m}}\bar{c}_{i}\Psi_i(\bx),
						\label{4.1:2} 
					\end{align}
					where 
					\[
					\bar{c}_{i}=\displaystyle\int_{\mathbb{A}^N}y_L(\bx)\Psi_i(\bx)f_{\bX}(\bx)\rd\bx
					\]
					and $\Psi_i(\bx)$, $i=1,\ldots,L_{N,S,m}$, are orthonormal polynomials consistent with the probability measure $f_{\bX}(\bx)\rd\bx$ of $\bX$. Then, the resulting $(\bar{m}+1)$-dimensional column vector of $\bM_{\bar{m}}(y_L)$ is 
					\begin{equation}	 
						\tilde{\bM}_{\bar{m}}(y_L)=(1,{\tilde{y}_{L,S,m}},\ldots,{\tilde{y}_{L,S,m}}^{\bar{m}})^{\intercal}.
						\label{4.1:3}
					\end{equation}	
					
					\item
					Construct an $(\bar{m}+1) \times (\bar{m}+1)$ monomial moment matrix of $\tilde{\bM}_{\bar{m}}({y}_L)$, defined as
					\begin{equation}
						\begin{split}
							\tilde{\bG}_{\bar{m}}:= \Exp[\tilde{\bM}_{\bar{m}}(Y_L)\tilde{\bM}_{\bar{m}}^{\intercal}(Y_L)]=
							\int_{{\bar{\mathbb{A}}}}\tilde{\bM}_{\bar{m}}(y_L)\tilde{\bM}_{\bar{m}}^{\intercal}(y_L)
							f_{Y_L}(y_L)\rd y_L.
						\end{split}
						\label{4.1:4}
					\end{equation}		
					For an arbitrary PDF $f_{Y_L}(y_L)$, $\tilde{\bG}_{\bar{m}}$ can be estimated with good accuracy by sampling methods, such as MCS, QMCS, or LHS, etc. 
					
					\item
					Select the $(\bar{m}+1) \times (\bar{m}+1)$ whitening matrix $\tilde{\bW}_{\bar{m}}$ from the Cholesky decomposition of the monomial moment matrix $\tilde{\bG}_{\bar{m}}$~\citep*{rahman2018polynomial}, leading to
					\begin{equation}
						\tilde{\bW}_{\bar{m}}^{-1}\tilde{\bW}_{\bar{m}}^{-\intercal}=\tilde{\bG}_{\bar{m}}.
						\label{4.1:5}
					\end{equation}	
					Then employ the whitening transformation to generate orthonormal polynomials from
					\begin{equation}
						\tilde{\mathbf{\Psi}}_{\bar{m}}(y_L)=\tilde{\bW}_{\bar{m}} \tilde{\bM}_{\bar{m}}(y_L).
					\end{equation}	
				\end{enumerate}
				
				These three steps are similar to those used for creating measure-consistent orthonormal polynomials in Section~\ref{sec:2.4:1}. However, in Steps 1 and 2 of the latter for $Y_L$, the realizations of $Y_L$ are determined by $S$-variate, $m$th-order DD-GPCE approximations. In doing so, the same basis functions $\Psi_i(\bX),~i=1,\ldots,L_{N,S,m}$, created in the three-step process in Section~\ref{sec:2.4:1} are reused. Furthermore, their respective DD-GPCE coefficients $\bar{c}_i,~i=1,\ldots,L_{N,S,m}$, are computed via SLS using the computationally economical lower-fidelity model. Therefore, the latter three-step process can be performed efficiently.  
				As discussed earlier for the related three-step process in Section \ref{sec:2.4:1}, constructing a well-conditioned version of the monomial moment matrix is critical in implementing the Cholesky factorization in \eqref{4.1:5}.

				\subsubsection{Expansion coefficients} \label{sec:4.2:2} 
				For $h:\mathbb{\bar{A}}\rightarrow \real$ and a realization $\bx=(x_1,\ldots,x_N)\in\mathbb{A}^N$ of $\bX$, let the input-output data set be 
				$$
				\{y_L^{(l)},y_H^{(l)}\}_{l=1}^{L^{\prime\prime}}:= \{y_L(\bx^{(l)}),y_H(\bx^{(l)})\}_{l=1}^{L^{\prime\prime}},
				$$
				with sample size $L^{\prime\prime}\in\nat$.  Then, the coefficients of the Fourier polynomial expansion in \eqref{4.1:1} are obtained by minimizing the mean-squared residual 
				\begin{align}
					e_{S,m}^{\prime\prime} :=
					\displaystyle
					\frac{1}{L^{\prime\prime}}
					\displaystyle
					\sum_{l=1}^{L^{\prime\prime}}
					\left[
					y_H^{(l)} -
					\displaystyle
					\sum_{i=1}^{\bar{m}+1} b_i \Psi_i(y_L^{(l)})
					\right]^2
					\label{4.1:6}
				\end{align}  
				via SLS explained in Section \ref{sec:2.4:3}. 
				
				\begin{remark}
					The Fourier-polynomial expansion in \eqref{4.1:6} includes scalar realizations $y_L^{(l)}$ of the random variable $Y_L$. The degree $\bar{m}$ is usually set as a small number, say $\bar{m}=1\text{--}3$. Thereby, the requisite sample size $L''$, usually determined ad hoc by requiring that $3\leq L''/(\bar{m}+1)\leq 8$, is small as well. This makes the bi-fidelity method to approximate the high-fidelity output data using $\bar{m}$th-order Fourier polynomial expansion computationally economical.  
				\end{remark}  
				
				\subsection{Complete algorithm and cost for $\mathrm{CVaR}$ estimation} \label{sec:4.3} 
				%
				Algorithm \ref{al2} presents a complete algorithm for the bi-fidelity method to estimate $\mathrm{VaR}_{\beta}$ and $\mathrm{CVaR}_{\beta}$.
				The cost to determine the $S$-variate, $m$th-order DD-GPCE approximation is dominated by the cost of evaluating the input-output data set $\{\bx^{(l)},y_H(\bx^{(l)})\}_{l=1}^{L^{\prime}}$ which requires $L^{\prime}\in\nat$ high-fidelity model evaluations. 
				In this regard, for high-dimensional problems, the DD-GPCE method via SLS alleviates the curse of dimensionality due to the sample size $L^{\prime}$ being determined as a multiple of the number of basis functions $L_{N,S,m}$ in \eqref{sr2} instead of the sample size of regular GPCE, $L_{N,m}$, in \eqref{2.3:1}. 
				Additionally, the proposed bi-fidelity method uses low-fidelity models to generate output data that approximate the high-fidelity model well, yielding an overall efficient method to compute the $S$-variate, $m$th-order DD-GPCE. 
				
				To select the fidelity for both models, the high-fidelity model is usually determined by a convergence test. A strategy for the selection of the low-fidelity models can be devised based on computational cost while maintaining physicality of the solution. Let $c_T$ be the total computational budget and $c_H, c_L$ be the costs of the high-fidelity and low-fidelity model evaluations, respectively.  Since the low and high-fidelity output sample sizes $L'$ and $L''$ are assumed to be inputs in Algorithm~\ref{al2}, the level of the low-fidelity model can be determined by considering its permissible cost  $c_L'$, i.e.,  
						\[
						c_L' \leq c_L= \dfrac{c_T-L''c_H}{L'}. 
						\label{cost} 
						\]
						Other strategies for selecting the low-fidelity model are discussed in Section~\ref{sec:5.2}. 
				\begin{algorithm} 
					\caption{Bi-fidelity method for $\mathrm{VaR}_{\beta}$ and $\mathrm{CVaR}_{\beta}$ estimation} \label{al2}
						\begin{algorithmic}[1]
							\Require Set a risk level $\beta \in (0,1)$. Set truncation parameters $S$, $m$, and $\bar{m}$. Set sample sizes $L$, $\bar{L}$, $L^{\prime}$, and $L^{\prime\prime}$. Generate input samples $\{\bx^{(l)}\}_{l=1}^L$ from the known probability measure $f_{\bX}(\bx)\rd\bx$ via MCS, QMCS, or LHS. 
							\Ensure  Estimate ${\mathrm{VaR}}_{\beta}[y(\bX)]$ and ${\mathrm{CVaR}}_{\beta}[y(\bX)]$ in the following steps: 
							\Procedure{Derive orthonormal (ON) polynomials}{$\mathbf{\Psi}_{S,m}(\bx)$}
							\State Set monomial vector $\bM_{S,m}(\bx)$ \eqref{2.4:1}.
							\State Create monomial moment matrix $\bG_{S,m}$ \eqref{2.4:2}.
							\State Create ON polynomial vector $\mathbf{\Psi}_{S,m}(\bx)$ by whitening transformation \eqref{2.4:3}.  
							\EndProcedure
							\Procedure{Obtain coefficients}{$\bar{c}_{i}$}
							\State Generate input-output data set $\{\bx^{(l)},y_L(\bx^{(l)})\}_{l=1}^{L^{\prime}}$ of size $L^{\prime}\in\nat$, where $y_L(\bx)$ presents a low-fidelity output.
							\State Use SLS to estimate $\bar{c}_{i},~i=1,\ldots,L_{N,S,m}$.  
							\EndProcedure
							\State Construct the $S$-variate, $m$th-order DD-GPCE $\tilde{y}_{L,S,m}(\bx)$ of low-fidelity output $y_L(\bx)$ \eqref{4.1:2}.
							\Procedure{Derive ON polynomials}{$\mathbf{\Psi}_{\bar{m}}(y_L)$}
							\State Set monomial vector $\tilde{\bM}_{\bar{m}}(y_L)$ via $\tilde{y}_{L,S,m}$  \eqref{4.1:3}.
							\State Construct monomial moment matrix $\tilde{\bG}_{\bar{m}}$ \eqref{4.1:4}.
							\State Create ON polynomial vector $\tilde{\mathbf{\Psi}}_{\bar{m}}(y_L)$ by whitening transformation \eqref{4.1:5}. 
							\EndProcedure
							\Procedure{Obtain coefficients}{$b_{i}$}
							\State Generate input-output data set $\{y_L^{(l)},y_H^{(l)}\}_{l=1}^{L^{\prime\prime}}$ of size $L^{\prime\prime}\in\nat$
							\State Use SLS to estimate $b_{i},~i=1,\ldots,\bar{m}+1$.  
							\EndProcedure
							\State Construct the $\bar{m}$th-order Fourier polynomial expansion $y_{H,\bar{m}}(\bx)$ of high-fidelity output $y_H(\bx)$ \eqref{4.1:1}. 
							\State By replacing $y$ with $y_{H,\bar{m}}$, construct $y_{S,m}$ \eqref{dd-gpce} and then perform Algorithm \ref{al1} to estimate $\mathrm{VaR}_{\beta}$ and $\mathrm{CVaR}_{\beta}$.
						\end{algorithmic}
				\end{algorithm}
				\section{Numerical results} \label{sec:5} 
				
				Two numerical examples are presented to illustrate the proposed DD-GPCE and bi-fidelity methods for estimating CVaR. In Section~\ref{sec:5.1}, a three-dimensional $36$-bar truss structure is considered, followed by an example of a glass/vinylester composite plate in Section~\ref{sec:5.2}. 

				In Examples~1~and~2, the sample size $L$ for sampling-based $\rm{CVaR}$ estimations is 10,000. This number is determined by a convergence test that yielded less than $0.1\%$ difference between the previous and current steps in MCS solutions for both examples. The monomial moment matrices $\bG_{S,m}$ and $\bG_{\bar{m}}$ in \eqref{2.4:2} and \eqref{4.1:4}, respectively, are determined by QMCS with $\bar{L}=5\times10^{6}$ samples together with the Sobol sequence~\citep*{sobol67}. The selection of the Sobol sequence in this work is due to its simplicity and efficiency in generating low-discrepancy quasi-random samples, thus improving the performance of QMCS. 

				In both examples, the coefficients of DD-GPCE or regular GPCE are estimated by SLS, c.f. Section~\ref{sec:2.4:3}. For a satisfactory estimation of the coefficients, we select the factors $L'/L_{N, S, m}$ (in DD-GPCE) or $L'/L_{N,m}$ (in regular GPCE) equal to three and four in Examples 1 and 2, respectively. In Example~2, we determine the number $L''$ of high-fidelity FEA by selecting the factor $L''/(\bar{m}+1)$ as eight.    
				The numerical results are obtained using MATLAB~\citep*{Matlab21} on an Intel Core i7-10850H 2.70 GHz processor with 64 GB of RAM.

				The $\rm{CVaR}$ estimates via the $S$-variate, $m$th-order DD-GPCE of $y(\bX)$ are denoted $\widehat{\mathrm{CVaR}}_{\beta}[y_{S,m}(\bX)]$. The proposed $\rm{CVaR}$ solutions are compared with a reference one by crude MCS of the chosen $L$ high-fidelity output data in a single trial.
				
				To measure the deviation of the proposed CVaR solution from that crude MCS estimate, we provide the mean relative difference (MRD) with respect to the crude MCS $\widehat{\mathrm{CVaR}}_{\beta}[y(\bX)]$, i.e., 
						\begin{align}
							\text{MRD}=
							\frac{\dfrac{1}{K}\displaystyle\sum_{k=1}^K\bigg|\widehat{\mathrm{CVaR}}_{\beta}[y(\bX)]-\widehat{\mathrm{CVaR}}_{\beta}[y_{S,m}^{(k)}(\bX)]\bigg|}{\bigg|\widehat{\mathrm{CVaR}}_{\beta}[y(\bX)]\bigg|},
							\label{5} 
						\end{align}
						where $\widehat{\mathrm{CVaR}}_{\beta}[y_{S,m}^{(k)}(\bX)]$ is the estimate obtained on the $k$th trial, and $K$ (the number of trials) is 20 in both Examples~1~and~2. At each trial, whether employing either the DD-GPCE, the regular GPCE, or Fourier-polynomial approximation for calculating DD-GPCE, we randomly select a subset of the $10,000$ high-fidelity output samples that are already available from the crude MCS. The MRD presents only the mean deviation of a set of proposed $\mathrm{CVaR}$ solutions over $K$ trials from a benchmark one obtained by crude MCS, which is not exact but approximate.
						We also report the average CVaR estimates over $20$ trials on Tables~\ref{table1} and \ref{table3} in Examples~1 and 2, respectively.
				\subsection {Example~1: A 36-bar 3D truss structure} \label{sec:5.1}
				This example demonstrates the efficacy of the DD-GPCE method in estimating the conditional value-at-risk of a system with a relatively high number ($N=36$) of dependent input random variables.   
				\subsubsection {Problem description}
				\begin{figure*}
					\begin{center}
						\includegraphics[angle=0,scale=0.65,clip]{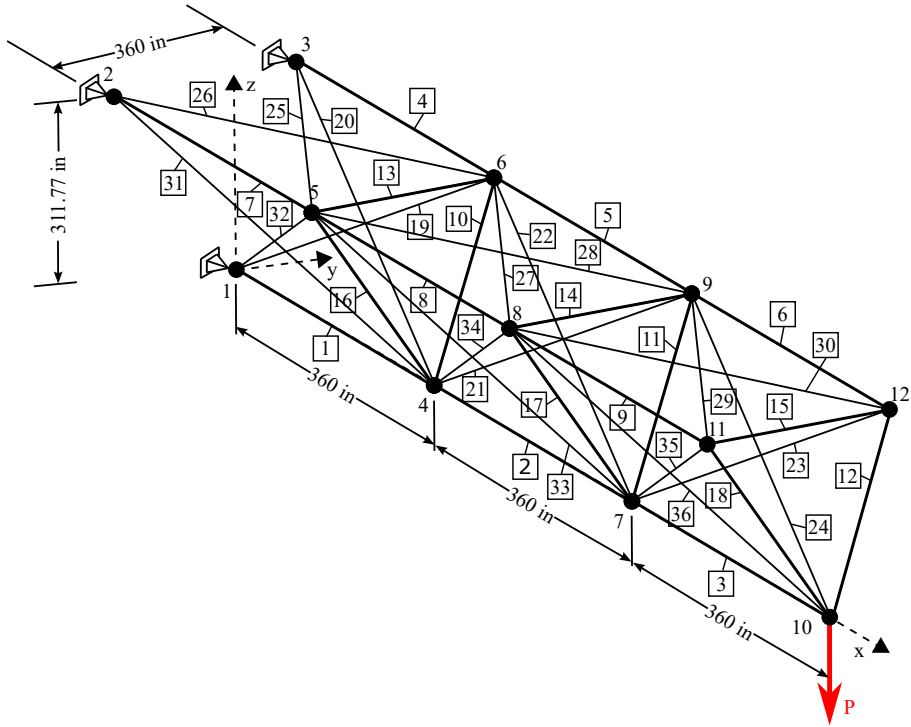}
					\end{center}
					\caption{Geometry, loading, and boundary conditions of the 36-bar 3D truss.}
					\label{fig1}
				\end{figure*}
				Figure~\ref{fig1} shows a 36-bar three dimensional truss that is simply supported at nodes 1, 2, and 3 and that is subject to a vertically downward concentrated force of 100 lb at node 10.  The truss is made of an aluminium alloy characterized by a Young's modulus of $10^7$ psi and a mass density of 0.1 lb/$\mathrm{in}^{3}$.  There are $N=36$ random variables $\bX=(X_1,\ldots,X_{36})^\intercal$ representing the respective random cross-sectional areas of thirty-six bars, that is, $X_i$ is the $i$th cross-sectional area, for $i=1,\ldots,36$. Described as multivariate Gaussian random variables, they have means $\Exp[X_i]=30$ and standard deviations equal to $0.05\Exp[X_i]$, $i=1,\ldots,36$; and correlation coefficients $\rho_{ij}=0.5$, $i,j=1,\ldots,36$, $i \ne j$. The problem is to determine $\mathrm{VaR}_{\beta}[y_l(\bX)]$ and $\mathrm{CVaR}_{\beta}[y_l(\bX)]$, $l=1,~2$, for two different quantile levels: $\beta=0.95$ (Case~1) and $\beta=0.99$ (Case~2). We consider two separate output quantities of interest, namely
				\begin{align}  
					y_1(\bX)&=\max\{ |u_1(\bX)|,\ldots,|u_{12}(\bX)|, |v_1(\bX)|,\ldots,|v_{12}(\bX)|\},
				\end{align} 			
				\begin{align}
					y_2(\bX)&=\max\{ |\sigma_1(\bX)|,\ldots,|\sigma_{36}(\bX)|\}.
				\end{align}
				Here, $|u_i|$ and $|v_i|$, are the absolute values of horizontal and vertical displacements, respectively, at the $i$th nodal point for $i=1\text{--}12$ and $|\sigma_{i}|$ is the absolute of axial stress of the $i$th bar for $i=1\text{--}36$. These quantities are determined via the spatially three-dimensional FEA for the linear elastic truss structure, which is implemented in Matlab with truss elements.
				The two cases of $\beta$ are shown to demonstrate the robustness of the DD-GPCE method for CVaR estimation with different risk levels.  
				\subsubsection {Results}
				\begin{table*}
					\caption{$\mathrm{CVaR}_{\beta}$ estimates  of two different output functions $y_1$ and $y_2$ of the 36-bar 3D truss structure via DD-GPCE and regular GPCE approximations and crude MCS in two distinct cases: Case 1 ($\beta=0.95$) and Case 2 ($\beta=0.99$): The sampling-based solutions $\widehat{\mathrm{CVaR}}_{\beta}$ are computed with a total of $L=10,000$ samples. }
					\begin{spacing}{1.2}
						\begin{center}
							\small
							\resizebox{\textwidth}{!}{
								\begin{tabular}{cccccccc}
									\toprule 
									\multirow{2}{*}{Method} &  & \multicolumn{2}{c}{Max. abs. displacement ($\mathrm{in}$) } &  & \multicolumn{2}{c}{Max. abs. axial stress ($\mathrm{psi}$) } & Number of FEA\tabularnewline
									\cmidrule{2-4} \cmidrule{3-4} \cmidrule{4-4} \cmidrule{6-7} \cmidrule{7-7} 
									&  & $\widehat{\mathrm{CVaR}}_{\beta}$ & MRD in \eqref{5}  &  & $\mathrm{CVaR}_{\beta}$ & MRD in \eqref{5} & for each trial \tabularnewline
									\midrule \tabularnewline 
									\multicolumn{8}{c}{Case 1 ($\beta=0.95$)}\tabularnewline
									1st-order regular GPCE$^{(\rm{a})}$ &  & $7.2102$$^{(\mathrm{b})}$ & $2.9533\times10^{-2}$ &  & $14010.1129$$^{(\mathrm{b})}$ & $3.8101\times10^{-2}$ & $111$\tabularnewline
									Univariate, 2nd-order DD-GPCE$^{(\rm{c})}$ &  & $7.2425$$^{(\mathrm{b})}$ & $2.5188\times10^{-2}$ &  & $14183.8278$$^{(\mathrm{b})}$ & $2.6174\times10^{-2}$ & $219$\tabularnewline
									Univariate, 3rd-order DD-GPCE$^{(\rm{d})}$ &  & $7.2439$$^{(\mathrm{b})}$ & $2.4992\times10^{-2}$ &  & $14240.7163$$^{(\mathrm{b})}$ & $2.2269\times10^{-2}$ & $327$\tabularnewline
									2nd-order regular GPCE$^{(\rm{e})}$ &  & $7.2435$$^{(\mathrm{b})}$ & $2.5049\times10^{-2}$ &  & $14153.3969$$^{(\mathrm{b})}$ & $2.8264\times10^{-2}$ & $2,109$\tabularnewline
									crude MCS &  & $7.2448$$^{(\mathrm{f})}$ & $-$ &  & $14157.2911$$^{(\mathrm{f})}$ & $-$ & $10,000$\tabularnewline
									\hdashline \tabularnewline
									\multicolumn{8}{c}{Case 2 ($\beta=0.99$)}\tabularnewline
									1st-order regular GPCE$^{(\rm{a})}$ &  & $7.3642$$^{(\rm{b})}$ & $8.7960\times10^{-3}$ &  & $14313.9756$$^{(\rm{b})}$  & $1.7239\times10^{-2}$ & $111$\tabularnewline
									Univariate, 2nd-order DD-GPCE$^{(\rm{c})}$ &  & $7.4233$$^{(\rm{b})}$  & $8.5367\times10^{-4}$ &  & $14633.5749$$^{(\rm{b})}$  & $6.2756\times10^{-3}$ & $219$\tabularnewline
									Univariate, 3rd-order DD-GPCE$^{(\rm{d})}$ &  & $7.4276$$^{(\rm{b})}$  & $2.7018\times10^{-4}$ &  & $14769.4983$$^{(\rm{b})}$  & $1.4554\times10^{-2}$ & $327$\tabularnewline
									2nd-order regular GPCE$^{(\rm{e})}$ &  & $7.4251$$^{(\rm{b})}$  & $6.0558\times10^{-4}$ &  & $14582.7015$$^{(\rm{b})}$  & $1.3539\times10^{-3}$ & $2,109$\tabularnewline
									crude MCS &  & $7.4296$$^{(\rm{f})}$  & $-$ &  & $14565.0584$$^{(\rm{f})}$  & $-$ & $10,000$\tabularnewline
									\bottomrule
							\end{tabular}}
							\par\end{center}
					\end{spacing}
					\begin{tablenotes}
						\scriptsize\smallskip 
						\item{a.} The first-order ($m=1$) regular GPCE is the same as the univariate ($S=1$), first-order ($m=1$) DD-GPCE.
						\item{b.} The ${\rm{CVaR}}$ estimate is averaged over 20 trials.  
						\item{c.} The truncation parameters are $S=1$ and $m=2$.	
						\item{d.} The truncation parameters are $S=1$ and $m=3$.				
						\item{e.} The second-order ($m=2$) regular GPCE is the same as the bivariate ($S=2$), second-order ($m=2$) DD-GPCE. 	
						\item{f.} The ${\rm{CVaR}}$ estimate is computed by crude MCS in one trial.
					\end{tablenotes}
					\label{table1}
				\end{table*}
				Table~\ref{table1} summarizes the sampling-based solution $\widehat{\mathrm{CVaR}}$ of $y_1(\bX)$ and $y_2(\bX)$ in Cases 1 and 2, including the requisite numbers of FEA by the regular GPCE and DD-GPCE methods. For comparison between the two methods, we also provide benchmark solutions in the form of crude MCS with 10,000 FEA as tabulated in the eighth and first rows from the bottom in Table~\ref{table1} in Cases~1 and~2, respectively.  
				The regular DD-GPCE and univariate DD-GPCE methods presented in Table~\ref{table1} yield estimates of CVaR of $y_1$ and $y_2$ that are very close to the crude MCS. Indeed, the range of their maximum value of MRD over $20$ trials is from $2.95\%$ to $3.81\%$ for the outputs of interest $y_1$ and $y_2$ for two cases of $\beta$. As expected, the first-order ($m=1$) regular GPCE approximations of $y_1$ and $y_2$ in both Cases 1 and 2 all show relatively lower precision than the other methods. 
				Observe that employing a second-order ($m=2$) regular GPCE approximation results in more precise CVaR estimates but the required FEA evaluations increase exponentially from 111 ($m=1$) to 2,109 ($m=2$). Therefore, for this high-dimensional ($N=36$) problem, the regular GPCE's curse of dimensionality becomes apparent. 
				On the other hand, the univariate ($S=1$) DD-GPCE approximations demand only 219--327 FEA evaluations as the degree $m$ increases from two to three; meanwhile, the accuracy of these solutions is better than those of the first-order ($m=1$) regular GPCE approximations. In particular, for $y_1$, the CVaR solutions by all univariate ($S=1$) DD-GPCE's are almost identical to those by the second-order ($m=2$) regular GPCE, but the DD-GPCE-based solutions provide a 5--10$\times$ cost savings compared to the second-order ($m=2$) regular GPCE. This illustrates the benefits of the proposed DD-GPCE method over the regular GPCE method in terms of computational efficiency for CVaR estimation with relatively high-dimensional dependent random input variables.  
				%
				\subsection {Example~2: A glass/vinylester composite plate} \label{sec:5.2}
				The second example focuses on an additional computational challenge, namely a nonlinear quasi-static FEA with a high number of internal states, where surrogate modeling is required to make the problem tracktable. Therefore, the proposed bi-fidelity method from Section~\ref{sec:4} is demonstrated. The number of dependent input random variables is moderate to high at $N=28$.   
				%
				\subsubsection{Problem description}
				\begin{figure*}
					\begin{center}
						\includegraphics[angle=0,scale=0.9,clip]{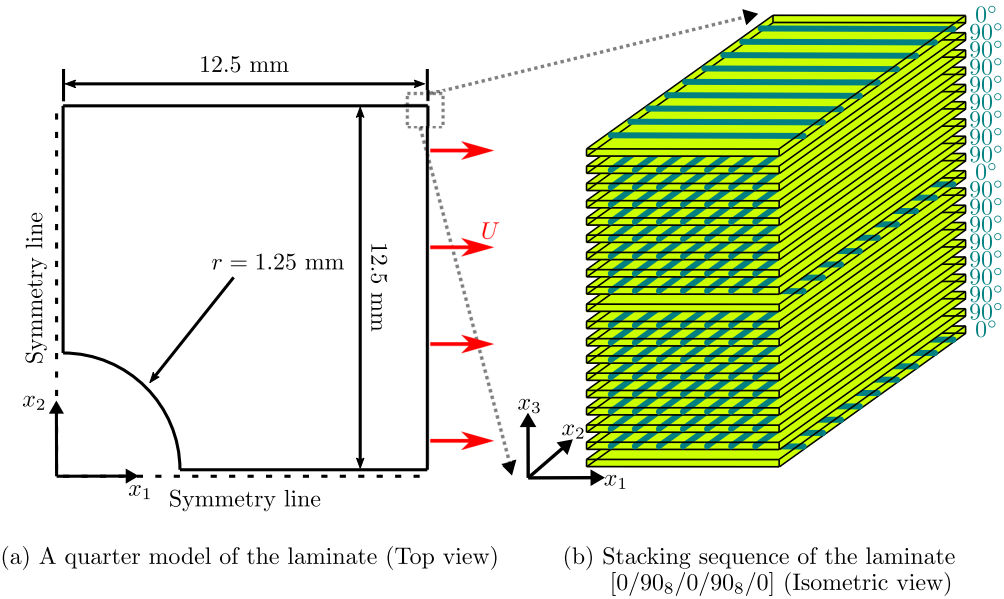}
					\end{center}
					\caption{Geometry, loading, and boundary conditions of a glass/vinylester composite plate: We consider one quarter of the full plate domain, which can represent the original square plate with domain $\mathcal{D}=(25~\mathrm{mm}\times25~\mathrm{mm})$ including a circular hole of radius $r=1.25~\mathrm{mm}$ in the center due to the symmetry conditions in the $x_1$ and $x_2$ directions of the center lines.}
					\label{fig2}
				\end{figure*}		
				Consider a two-dimensional square plate with spatial domain $\mathcal{D}=(25~\mathrm{mm}\times25~\mathrm{mm})$, including a circular hole of radius $r=1.25~\mathrm{mm}$ in the center. Due to symmetry, we consider only a quarter of the plate, resulting in the domain $\bar{\cD}=(12.5~\mathrm{mm}\times12.5~\mathrm{mm})$, as shown in Figure~\ref{fig2}a. 
				Figure~\ref{fig2}b shows that the glass/vinylester laminate (Fiberite/HyE 9082Af) is composed of 19 stacking sequences $[0/90_8/0/90_8/0]$, where `0' indicates a ply having fibers in $x_1$ direction and $90_8$ indicates that eight plies have fibers in $x_2$ direction.
				The plate is subjected to a uniaxial tensile displacement loading $U$ that acts on the entire right side edge. The geometry of the composite plate was initially studied by~\cite*{MOURE2014224}. 
				
				The local directions are $x_1,~x_2,~\text{and}~x_3$ and we define $E_i$, $\nu_{ij}$, and $G_{ij}$ to be Young's modulus, Poisson's ratio, and shear modulus of the plies in the corresponding local directions for $i,j=1,2,3$. Also, let $S_{ti}$, $S_{ci}$, and $S_{sij}$ be tensile, compressive, and shear strengths of the plies. The measured material properties and ply thicknesses of the laminate vary, and we therefore model them as $N=28$ random variables, as presented in Table~\ref{table2}. The $19$ random variables that model the ply thickness are modeled as correlated via a multivariate lognormal distribution with a correlation coefficient of $0.5$. The remaining $9$ random variables are modeled as independent and distributed uniformly.     
				\begin{table*}
					\caption{Statistical properties of constituents in Fiberite/HyE 9082Af}
					\begin{spacing}{1.0}
						\begin{centering}
							\footnotesize
								\begin{tabular}{ccccccc}
									\toprule 
									Random  & \multirow{2}{*}{Property} & \multirow{2}{*}{Mean} & \multirow{2}{*}{COV (\%)} & Lower  & Upper & Probability \tabularnewline
									variable &  &  &  & boundary & boundary & distribution\tabularnewline
									\midrule 
									$X_{1}$ & $E_{1}$ (MPa) & $44,700$ & 11.55  & 35760 & 53640 & Uniform\tabularnewline
									$X_{2}$ & $E_{2}$ (MPa) & $12,700$ & 11.55  & 10,160 & 15,240 & Uniform\tabularnewline
									$X_{3}$ & $v_{12}$ & $0.297$ & 11.55  & 0.238 & 0.356 & Uniform\tabularnewline
									$X_{4}$ & $G_{12}$ (MPa) & $5,800$ & 11.55  & 4,640 & 6,960 & Uniform\tabularnewline
									$X_{5}$ & $S_{t1}$ (MPa) & $1,020$ & 11.55  & 816 & 1,224 & Uniform\tabularnewline
									$X_{6}$ & $S_{t2}$ (MPa) & $40$ & 11.55  & 32 & 48 & Uniform\tabularnewline
									$X_{7}$ & $S_{c1}$ (MPa) & $620$ & 11.55  & 496 & 744 & Uniform\tabularnewline
									$X_{8}$ & $S_{c2}$ (MPa) & $140$ & 11.55  & 126 & 168 & Uniform\tabularnewline
									$X_{9}$ & $S_{s12}$$^{(b)}$ (MPa) & $60$ & 11.55  & 48 & 72 & Uniform\tabularnewline
									\multirow{2}{*}{$X_{10}$--$X_{28}$$^{(\mathrm{a})}$} & Plies 1--19 & \multirow{2}{*}{$0.144$} & \multirow{2}{*}{$6$} & \multirow{2}{*}{0} & \multirow{2}{*}{$\infty$} & Multivariate\tabularnewline
									& thicknesses &  &  &  &  & Lognormal\tabularnewline
									\bottomrule
								\end{tabular}
								\par\end{centering}
						\end{spacing}\tabularnewline
						\begin{tablenotes}
							\scriptsize\smallskip
							\item{a.} Correlation coefficients among $X_{10}$--$X_{28}$ are 0.5.
							\item{b.} $S_{s12}=S_{s23}.$
						\end{tablenotes}			
						\label{table2}
					\end{table*}
					%
					\subsubsection{Hashin damage criterion}
					To perform a quasi-static damage analysis of the fiber-reinforced composite laminate we choose the Hashin damage model. This captures the nonlinear behavior of the composite material during failure progression, which in turn allows for a detailed description of the failure mechanisms.  
					The damage model includes a damage initiation criterion that takes into account four possible failure modes~\citep*{DUARTE2017277}: (1) fiber breakage in tension ($F_f^t$), (2) fiber buckling in compression ($F_f^c$), (3) matrix cracking in tension ($F_m^t$), and (4) matrix crushing in compression ($F_m^c$).  A fiber-reinforced composite is considered damaged if 
					\begin{align}\label{5.2:1}
						F_f^t=\left(\dfrac{\sigma_{11}}{S_{t1}}\right)^2 + \alpha\left(\dfrac{\sigma_{12}}{S_{s12}}\right)^2 \geq 1.0~\text{and}~\sigma_{11}\ge 0,  
					\end{align}
					\begin{align} 
						F_f^c=\left(\dfrac{\sigma_{11}}{S_{c1}}\right)^2 \geq 1.0~\text{and}~\sigma_{11} < 0,  
					\end{align}
					\begin{align} 
						F_f^t=\left(\dfrac{\sigma_{22}}{S_{t2}}\right)^2 + \left(\dfrac{\sigma_{12}}{S_{s12}}\right)^2 \geq 1.0~\text{and}~\sigma_{22}\ge 0,  
					\end{align}
					and
					\begin{align}
						F_m^c=\left(\dfrac{\sigma_{22}}{2S_{s23}}\right)^2 + \left[\left(\dfrac{S_{c2}}{2S_{s23}}\right)^2 - 1 \right]\dfrac{\sigma_{22}}{S_{c2}} + \left(\dfrac{\sigma_{12}}{S_{s12}}\right)^2 
						\geq 1.0~\text{and}~\sigma_{22}\ge 0,
					\end{align}
					where $\sigma_{11}$ and $\sigma_{22}$ are the principal stresses in the $x_1$ and $x_2$ directions, respectively, and $\sigma_{12}$ is the shear stress in $x_1$--$x_2$ plane. Also, the coefficient $\alpha\in [0,1]$ in \eqref{5.2:1} accounts for the contribution of the shear stress $\sigma_{12}$ to the fiber breakage, and set as $\alpha=1$ in this work. 
					We employ the Hashin criterion together with the progressive damage model; both are built-in functions in ABAQUS/Explicit, version~6.14-2.
					%
\subsubsection{Output of interest}
\begin{figure*}
\begin{center}
\includegraphics[angle=0,scale=0.80,clip]{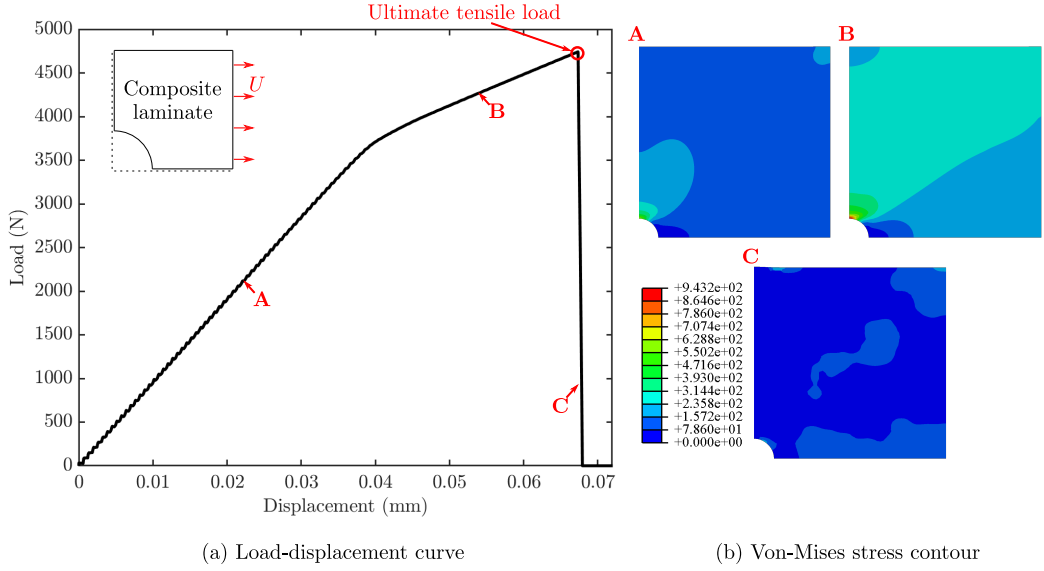}
\end{center}
\caption{FEA results for the glass/vinylester composite plate: The ultimate tensile load in the load-displacement curve (a) is recorded as the output of interest. Part (b) shows von-Mises stress contour obtained in regimes A, B, and C of the load displacement curve, indicating that the stress concentration occurred at the top of circular hole advances in the regime A to  B before a fracture occurs in the regime C.}
\label{fig3}
\end{figure*}	
The ultimate tensile load is chosen as the output quantity of interest.  Figure~\ref{fig3} describes the resulting load versus applied displacement curve of the composite laminate. As the applied displacement at the right edge increases in the range of $0$--$0.04$~$\mathrm{mm}$, the resulting load almost linearly increases, which we label as regime `A' in the figure.  However, after the displacement reaches $0.04~\mathrm{mm}$, the rate of change of the resulting load over the displacement decreases, which we indicate as regime `B' in the figure.  This is because the matrix or fiber starts to be degraded as a part of them exceeds the Hashin damage criteria. Thereafter, it reaches a peak value until it drops dramatically. This indicates complete fracture, labeled  as regime `C'.  Indeed, Figure~\ref{fig3}b  supports the above description of the damage at three marks in Figure~\ref{fig3}a. In regime A, the stress concentration appears at the top of the hole, and the stress is further advanced in regime B. Finally, complete damage occurs, fracturing the composite laminate in the $x_2$ direction from the top of the circular hole and losing the stress distribution at C, as shown in Figure~\ref{fig3}b.     
					
Additionally, to confirm the accuracy of the Hashin damage model, we compare the ultimate tensile load ($130.45~\mathrm{MPa}$) by the numerical model for the T300/1034-C laminate with stacking sequence $[0/\pm40/90_7/90_7/\mp40/0]$ with experimental results ($134.5~\mathrm{MPa}$) reported in the literature~\citep*{1991Tan}. The relative error between experiment and simulation is only 3\%. 					 
					
\subsubsection{High-fidelity and low-fidelity models}
\begin{figure*}
\begin{center}
\includegraphics[angle=0,scale=0.80,clip]{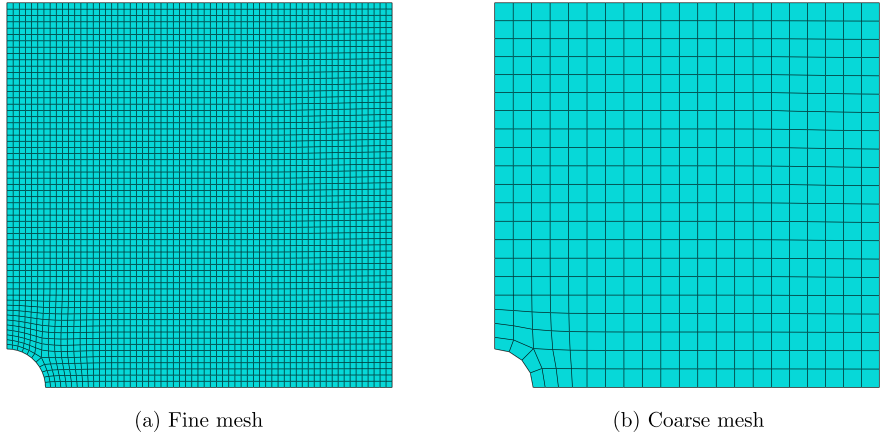}
\end{center}
\caption{The finite element mesh of the glass/vinylester composite plate: The fine mesh in (a) comprising 3,887 elements is used to generate high-fidelity output data, while the coarse mesh in (b) comprising 441 elements is used to generate lower-fidelity output data.}
\label{fig4}
\end{figure*}

To evaluate the bi-fidelity method presented in Section~\ref{sec:4}, the high and low-fidelity outputs were obtained by FEA as follows: the high-fidelity output is computed from a fine mesh model with 3,887 elements, and the lower-fidelity model from a coarse mesh with only 441 elements, as presented in Figures \ref{fig4}a and b, respectively. The element types were chosen as S4R, a 4-node, quadrilateral, stress/displacement shell element with reduced integration, built in ABAQUS/Explicit, version 6.14-2. Therefore, the total number of degrees of freedom for the high-fidelity and lower-fidelity models are 24,084 and 2,910, respectively. The global mesh size ($0.2$ mm) of the high-fidelity model was determined on the basis of a convergence study, i.e., refining the mesh size globally until the output quantity of interest becomes converged. On the other hand, we determined the low-fidelity model by selecting a larger global mesh size ($0.6$ mm), where the selected low-fidelity model must be able to capture an ultimate tensile load induced by damage (see Figure~\ref{fig3}) at given loading and boundary conditions. Alternatively, the coarse mesh size could be determined in a way that satisfies a given computational budget $c_T$ \eqref{cost}, while still maintaining that the solution is physical and realistic. Optimal load management in the spirit of~\citep*{peherstorfer2016optimal} for the surrogate model may also be adapted to CVaR computations.

\subsubsection{Results}
\begin{table*}
\caption{$\mathrm{CVaR}_{\beta}$ estimates ($\beta=0.99$) of the ultimate tensile load of the glass/vinylester composite plate via DD-GPCE approximations and the crude MCS: The sampling-based solutions $\widehat{\mathrm{CVaR}}_{\beta}$ are computed with a total of $L=10,000$ samples.}
\begin{spacing}{1.2}
\begin{center}
\small
\resizebox{\textwidth}{!}{
\begin{tabular}{clccccccc}
										\toprule 
										\multicolumn{3}{c}{} &  & \multicolumn{2}{c}{Ultimate tensile load (N)} & \multicolumn{2}{c}{Number of FEA in a single trial} & CPU time in\tabularnewline
										\multicolumn{3}{c}{Methods} &  & $\widehat{\mathrm{CVaR}}_{\beta}$ & MRD in \eqref{5} & High-fidelity$^{\rm{(a)}}$ & Low-fidelity$^{\rm{(b)}}$ &a single trial (hours)$^{\rm{(c)}}$\tabularnewline
										\midrule 
										\multicolumn{9}{l}{Univariate ($S=1$), third-order ($m=3$) DD-GPCE}\tabularnewline
										& \multicolumn{2}{l}{High-fidelity output} &  & $5540.6511$$^{\rm{(d)}}$ & $7.4339\times10^{-3}$ & $340$ & $-$ & $29.3$\tabularnewline
										& \multicolumn{2}{l}{Low-fidelity output} &  & $6141.2138$$^{\rm{(d)}}$ & $1.0015\times10^{-1}$ & $-$ & $340$ & $6.8$\tabularnewline
										& \multicolumn{2}{l}{Bi-fidelity approximation} &  &  &  &  &  & \tabularnewline
										&  & First-order ($\bar{m}=1$) Fourier-polynomial expansion &  & $5437.7325$$^{\rm{(d)}}$ & $2.7375\times10^{-2}$ & $16$ & $340$ & $8.1$\tabularnewline
										&  & Second-order ($\bar{m}=2$) Fourier-polynomial expansion &  & $5451.5180$$^{\rm{(d)}}$ & $2.9232\times10^{-2}$ & $24$ & $340$ & $8.8$\tabularnewline
										&  & Third-order ($\bar{m}=3$) Fourier-polynomial expansion &  & $5488.0055$$^{\rm{(d)}}$ & $2.1806\times10^{-2}$ & $32$ & $340$ & $9.5$\tabularnewline\tabularnewline\multicolumn{9}{l}{{Second-order ($m=2$) regular GPCE (equivalent to bivariate ($S=2$), second-order ($m=2$) DD-GPCE)}}\tabularnewline
										& \multicolumn{2}{l}{High-fidelity output} &  & $5576.9167^{\rm{(d)}}$ & $9.5682\times10^{-4}$ &{{ $1,740$ }}& $-$ & {{$151.6$}}\tabularnewline
										& \multicolumn{2}{l}{{Low-fidelity output}} &  & {{ $6157.3941^{\rm{(d)}}$}} & {{$1.0305\times10^{-1}$}} & $-$ & {{$1,740$}} & {{$36.1$}}\tabularnewline
										& \multicolumn{2}{l}{ {Bi-fidelity approximation}} &  &  &  &  &  & \tabularnewline
										&  &  {{First-order ($\bar{m}=1$) Fourier-polynomial expansion }}&  &  {{$5463.5153^{\rm{(d)}}$}} & {{$2.7798\times10^{-2}$}} &  {{$16$}} & {{$1,740$}} &  {{$37.5$}}\tabularnewline
										&  & {{Second-order ($\bar{m}=2$) Fourier-polynomial expansion}} &  &  {{$5442.6036^{\rm{(d)}}$}} &  {{$3.7529\times10^{-2}$}} & {{$24$}} & {{$1,740$}} & {{$38.2$}}\tabularnewline
										&  & {{Third-order ($\bar{m}=3$) Fourier-polynomial expansion}} &  &  {{$5472.6904^{\rm{(d)}}$}} &  {{$2.7351\times10^{-2}$}} &  {{$32$}} &  {{$1,740$}} &  {{$38.9$}}\tabularnewline\tabularnewline
										\multicolumn{3}{l}{crude MCS} &  & $5582.1484$$^{\rm{(e)}}$ & $-$ & $10,000$ & $-$ & $859.8$\tabularnewline
										\bottomrule
								\end{tabular}}
								\par\end{center}
\end{spacing}
\begin{tablenotes}
\scriptsize\smallskip
							\item{a.} The high-fidelity output is obtained by the fine mesh model in Figure~\ref{fig4}a.		
							\item{b.} The lower-fidelity output is obtained by the coarse mesh model in Figure~\ref{fig4}b.	
							\item{c.} The CPU time is the sum of the CPU time spent executing Algorithm~\ref{al2} including the number of FEA$\times$ FEA computation time in a single trial. The FEA computation time is averaged over five simulations. 
							\item{d.} The $\rm{CVaR}$ estimate is averaged over $K=20$ trials.  
							\item{e.} The $\rm{CVaR}$ estimate is computed by crude MCS in one trial.  	
\end{tablenotes}
\label{table3}
\end{table*}
					Table~\ref{table3} summarizes the CVaR estimates obtained from univariate ($S=1$) third-order ($m=3$) DD-GPCE and second-order ($m=2$) regular GPCE (which is equivalent to the bivariate ($S=2$) second-order ($m=2$) DD-GPCE) approximations. Each DD-GPCE, as its three distinct versions,  uses (i) the high-fidelity model, (ii) the lower-fidelity model and (iii) the bi-fidelity approximation based on Fourier-polynomial expansions, as presented in Section~\ref{sec:4.2}. 
					Compared to the standard MCS with 10,000 high-fidelity model evaluations (first row from the bottom of Table~\ref{table3}) the univariate ($S=1$) DD-GPCE-based CVaR estimates using the high-fidelity model yield $0.743\%$ in MRD. Notably, through proper pruning of basis functions, the univariate ($S=1$) DD-GPCE-based CVaR estimate requires only $340$ FEA (high-fidelity) solutions, leading to a speedup of $29.4\times$ compared to the $859.8~\mathrm{hours}$ of CPU time required for the standard MCS. 
					However, even that CPU time of $29.2~\mathrm{hours}$ for the univariate ($S=1$) DD-GPCE-based CVaR estimate via the high-fidelity model can still be computationally expensive when users want rapid turnaround in the design processes or when used for optimization problems.  
					We also report the univariate ($S=1$) DD-GPCE solution when only the lower-fidelity model is used. This requires a CPU time of 6.8 hours, almost one-fifth of those 29.3 hours required by the high-fidelity output but yields $10.015\%$ in MRD.  Thus, the univariate ($S=1$) DD-GPCE employing only the lower-fidelity output model yields an inaccurate or biased CVaR estimate, so it is of limited use. 
					
					To achieve \textit{almost} the accuracy of the high-fidelity CVaR solution but with the efficiency of the lower-fidelity version, we proposed and tested the proposed bi-fidelity method based on Fourier-polynomial expansion. 
					Indeed, as reported in the {eleventh through thirteenth} rows from the bottom of the Table~\ref{table3}, all CVaR estimates via the bi-fidelity method are very close to the benchmark CVaR estimate via the crude MCS. 
					As the degree ($\bar{m}$) increases, the precision of the respective CVaR solution increases from 2.738 \% to 2.181 \% in MRD, but it shows some fluctuation at the case of $\bar{m}=2$ due to a random sampling effect. The respective average CVaR solution (the eleventh through thirteenth rows from the bottom in the second column of Table~\ref{table3}) monotonically approaches the benchmark value.

					Remarkably, a $2.738~\%$ in MRD can be achieved while needing $8.0$ hours in CPU time, which is only marginally higher than the $6.7$ hours demanded by the low-fidelity version but provides an almost five-fold improvement in accuracy. Moreover, that third-order ($\bar{m}=3$) bi-fidelity approximation yields the precision of $2.181~\%$ in MRD by requiring only an additional $32$ high-fidelity FEA alongside many lower-fidelity evaluations. This situation is typical for multifidelity UQ methods, where the bulk of the computations are done by lower-fidelity approximations, yet a few high-fidelity solutions provide a significant boost of overall accuracy (and other benefits, see~\cite{peherstorfer2018survey}).
					Moreover, the second-order ($m=2$) regular GPCE (equivalent to the bivariate ($S=2$) second-order ($m=2$) DD-GPCE) approximations provide similar results compared to the univariate ($S=1$) third-order ($m=3$) version, as presented in the third through eighth rows of Table~\ref{table3}.  In other words, while the bivariate ($S=2$) DD-GPCE-based CVaR estimate using the high-fidelity model is the most accurate, showing $0.096\%$ in MRD, its counterpart using the low-fidelity model yields the lowest accurate solution ($10.31\%$ in MRD). However, the requisite number ($1,740$) of FEA by the bivariate ($S=2$), second-order ($m=2$) DD-GPCE is more than five times of those ($340$) by its univariate ($S=1$), third-order ($m=3$) version. As a result, the bivariate DD-GPCE based solution using the high-fidelity model demands almost 151.6 hours in CPU time in each trial which is also more than five times expensive than one (29.3 hours) by the univariate version using the high-fidelity model.  
					
					To reduce the 151.6 hours CPU time to the level ($36.1$ hours) of the low-fidelity version but obtain a more precise CVaR solution than the low-fidelity based solution (10.3\% in MRD), we employ the bi-fidelity method. While using first-order ($\bar{m}=1$) through third-order ($\bar{m}=3$) approximations of the bi-fidelity method, we obtained more accurate CVaR solutions ($2.7\%$ to $3.7\%$ in MRD) than the low-fidelity based solution ($10.31\%$ in MRD). Consequently, this result proves that the bi-fidelity method can be effective and robust to the case of bivariate ($S=2$) DD-GPCE or regular GPCE approximations.
					
We note that the CPU times (37.5--38.9 hours) by the bivariate ($S=2$), second-order ($m=2$) DD-GPCE using the bi-fidelity approximations are higher than those (29.3 hours) by the univariate ($S=1$), third-order ($m=3$) DD-GPCE using high-fidelity model. However, the accuracy of the former with bi-fidelity approximations is less than the latter with the high-fidelity model. This result indicates that when a higher value of $S$ or $m$ is chosen than needed to satisfy the user-defined threshold, the efficiency gains diminish for high-dimensional inputs. This is due to rapidly growing number $L_{N,S,m}$ of the coefficients or basis functions as $S$ or $m$ increases. Having said this, the bi-fidelity method improves the accuracy of the CVaR solution significantly (from 10\% to 2--3\% in MRD) while only requiring slightly more CPU resources (3.9\%--7.8\%) than the low-fidelity version. These results demonstrate both the robustness of the bi-fidelity method and the power of the DD-GPCE to improve its efficiency by truncating the basis functions in a dimension-wise manner using $S$. 							
\begin{figure*}
\begin{center}
\includegraphics[angle=0,scale=0.85,clip]{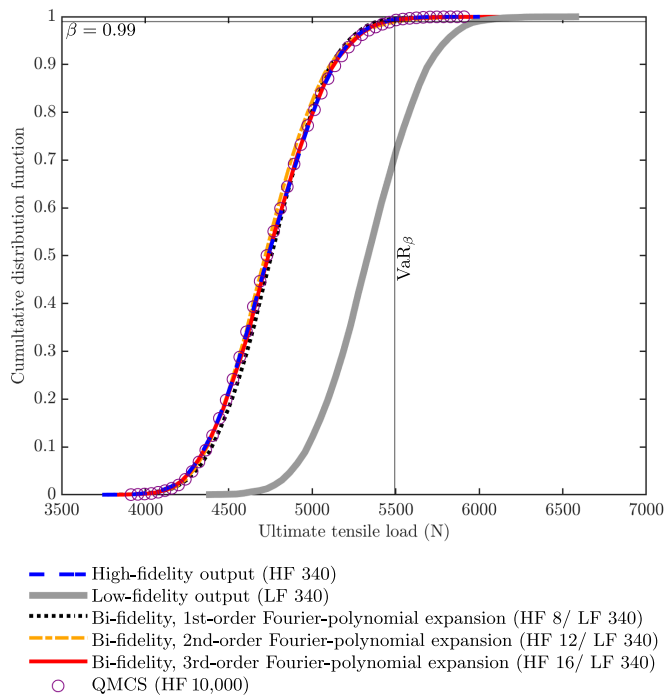}
\end{center}
\caption{Cumulative distribution function of the ultimate tensile load in the glass/vinylester composite plate: High-fidelity (HF) and low-fidelity (LF) output and three bi-fidelity approximations are used for univariate, $3$rd-order DD-GPCE approximations. The calculated DD-GPCE approximations are resampled 10,000 times to estimate the CDFs.}
\label{fig5}
\end{figure*}					
Figure~\ref{fig5} shows the estimated CDFs of the ultimate tensile load by three distinct versions of the univariate ($S=1$), third-order ($m=3$) DD-GPCE approximations, i.e., via high- and low-fidelity output data and also via the bi-fidelity approximation. For each of these distinct methods, a single trial CDF estimate was chosen among the CDF estimates from $K=20$ trials. For reference, we also show the crude MCS obtained from 10,000 high-fidelity FEA simulations. 
When the DD-GPCE approximation is computed by the high-fidelity output data of size $340$, its CDF is almost identical to the CDF of the crude MCS. On the other hand, the DD-GPCE (of size $340$) with only the low-fidelity output produces a CDF  that would make an engineer deduce that the composite can withstand much higher tensile loads, which is not true. This can have detrimental effects in engineering practice. 
In contrast, the proposed DD-GPCE methods in conjunction with the bi-fidelity approximation provide CDF estimates that are very close to those by the crude MCS but require only $16$--$32$ high-fidelity FEA evaluations paired with $340$ low-fidelity FEA.
Specifically, as the degree $\bar{m}$ of the Fourier-polynomial expansion increases, the respective CDF estimates approach those by the high-fidelity version or the crude MCS very closely. This demonstrates the ability of the proposed bi-fidelity method to accurately estimate CDFs as well.
\section{Conclusion and future direction}  \label{sec:6}  
We presented a novel bi-fidelity method for efficient CVaR estimation of complex nonlinear systems subject to arbitrary and high-dimensional dependent random input variables. The new method entails (1) the DD-GPCE approximation of a stochastic output function with high-dimensional dependent inputs, (2) an innovative method employing Fourier-polynomial expansions of a mapping between the stochastic lower-fidelity and high-fidelity output data for efficiently calculating the DD-GPCE, and (3) a standard sampling-based CVaR estimation integrated with the DD-GPCE. 
The proposed bi-fidelity method uses measure-consistent orthonormal polynomials in the random variable of the low-fidelity output to approximate the high-fidelity output, thus achieving a nearly exponential convergence rate for the output data.
The strength of the bi-fidelity approach is that it only requires a handful of high-fidelity output evaluations to augment the (many) lower-fidelity evaluations. When equipped with high-order ($\bar{m}$) basis functions, the Fourier-polynomial expansions can achieve more accurate high-fidelity output approximations. 
The numerical results for the truss structure with $N=36$ (dependent) random variables showed that CVaR estimation is possible by combining sampling-based CVaR estimation with DD-GPCE approximations. The example demonstrates that DD-GPCE via a dimension-wise reconstruction of the GPCE basis functions can alleviate the curse of dimensionality over the regular GPCE when a high-dimensional problem is addressed. Finally, the power of the proposed bi-fidelity method to achieve nearly the accuracy of the high-fidelity CVaR solution with the efficiency of the low-fidelity version was demonstrated by solving the glass/vinylester laminate problem involving 28 (mostly dependent) input random variables.  
					
The DD-GPCE surrogate can be somewhat limited in terms of accuracy when high-variate interaction effects among input variables are not negligible to the output quantity of interest (i.e., those terms in the PCE cannot be pruned). While the proposed bi-fidelity method uses Fourier-polynomial expansions to approximate the high-fidelity output, it is agnostic to the surrogate model structure. Therefore, in cases where high-variate interactions are relevant, one may need to use and/or develop new surrogate modeling techniques to approximate the output. 
A potential for improvement of the proposed method is related to the selection of the truncation parameters, e.g., $S$ and $m$ for DD-GPCE or $\bar{m}$ for Fourier-polynomial expansion. Instead of choosing the truncation parameters arbitrarily, one may exploit an adaptive version of DD-GPCE or Fourier-polynomial expansion, where a truncated set of basis is chosen optimally based on a specified error tolerated by the resulting approximation.
Moreover, while the bi-fidelity method in this work was mainly used to estimate the value at-risk and the conditional value-at-risk, it can be extended to a more general class of uncertainty quantification problems, such as second-moment or reliability analysis.  
					
\section*{Acknowledgment} \label{sec:7} 

\noindent\textbf{\textit{Funding}} This material is based on research sponsored by the Air Force Research Lab (AFRL) and Defense Advanced Research Projects Agency (DARPA) under agreement number FA8650-21-2-7126. The U.S. Government is authorized to reproduce and distribute reprints for Governmental purposes notwithstanding any copyright notation thereon.
\newline\newline 
\noindent\textbf{\textit{Disclaimer}} The views and conclusions contained herein are those of the authors and should not be interpreted as necessarily representing the official policies or endorsements, either expressed or implied, of the AFRL and DARPA or the U.S. Government.  
					
\bibliography{ref}
\bibliographystyle{abbrv}
						
\appendix 
\renewcommand{\thesection}{{Appendix}~\Alph{section}}
	\section{Three step process to construct measure-consistent orthonormal polynomials}  \label{appendix:a} 
	This appendix summarizes a process to generate the multivariate orthonormal polynomial basis of GPCE in Section \ref{sec:2.3}. The orthonormal polynomial functions are consistent with an arbitrary, non-product-type probability measure $f_{\bX}(\bx)\rd\bx$ of $\bX$ and determined by the following three steps. 
	\begin{enumerate}
		\item
		Given $m \in \nat_0$, create an $L_{N,m}$-dimensional column vector
		\begin{equation}	 \bM_m(\bx)=(\bx^{\bj^{(1)}},\ldots,\bx^{\bj^{(L_{N,m})}})^{\intercal},
		\end{equation}
		of monomials whose elements are the monomials $\bx^{\bj}$ for $|\bj|\le m$ arranged in the aforementioned order.  It is referred to as the monomial vector in $\bx=(x_1,\ldots,x_N)^\intercal$ of degree at most $m$.
		
		\item
		Construct an $L_{N,m} \times L_{N,m}$ monomial moment matrix of $\bM_m(\bX)$, defined as
		\begin{equation}
			\begin{split}
				{\bG}_m:= \Exp[\bM_m(\bX)\bM_{m}^{\intercal}(\bX)]=
				\int_{\mathbb{\bar{A}}^{N}}\bM_{m}(\bx)\bM_{m}^{\intercal}(\bx)
				f_{\bX}(\bx)\rd\bx.
			\end{split}
		\end{equation}		
		For an arbitrary PDF $f_{\bX}(\bx)$, $\bG_m$ cannot be determined exactly, but it can be estimated with good accuracy by numerical integration and/or sampling methods \citep*{lee2020practical}.
		
		\item
		Select the $L_{N,m} \times L_{N,m}$ whitening matrix ${\bW}_m$ from the Cholesky decomposition of the monomial moment matrix ${\bG}_m$ \citep*{rahman2018polynomial}, leading to
		\begin{equation}
			{\bW}_{m}^{-1}{\bW}_{m}^{-\intercal}={\bG}_{m}.
		\end{equation}	
		Then employ the whitening transformation to generate multivariate orthonormal polynomials from
		\begin{equation}
			{\mathbf{\Psi}}_m(\bx)={\bW}_m \bM_m(\bx).
		\end{equation}	
	\end{enumerate}
\end{document}